\documentclass[opre,nonblindrev]{informs3} 

\OneAndAHalfSpacedXI 

\usepackage{amsmath}
\usepackage{amsfonts}
\usepackage{amssymb}
\usepackage{mathtools}
\usepackage{dsfont}
\usepackage{color}
\usepackage{multirow}
\usepackage[table,xcdraw]{xcolor}
\usepackage{float}
\usepackage{placeins}
\usepackage{afterpage}
\usepackage{caption}
\usepackage{fullpage}
\usepackage{epsfig}
\usepackage{algcompatible}
\usepackage{algorithm}
\usepackage{url}
\usepackage{pifont}
\usepackage{longtable}
\usepackage{tabularx}
\usepackage{pdflscape}
\usepackage{booktabs}

\newcommand{\ignore}[1]{}
\def\defi{\vcentcolon=}
\renewcommand{\theequation}{\thesection.\arabic{equation}}
\newcolumntype{P}[1]{>{\arraybackslash}p{#1}}

\usepackage{natbib}
 \bibpunct[, ]{(}{)}{,}{a}{}{,}%
 \def\bibfont{\small}%
\TheoremsNumberedThrough     
\ECRepeatTheorems

\EquationsNumberedThrough    

\allowdisplaybreaks

\begin{document}


\RUNAUTHOR{Guigues, Kleywegt, and Nascimento}

\RUNTITLE{Operation of an Ambulance Fleet}

\TITLE{Optimization-based Operation of an Ambulance Fleet}

\ARTICLEAUTHORS{%
\AUTHOR{Vincent Guigues}
\AFF{School of Applied Mathematics, FGV, Praia de Botafogo, Rio de Janeiro, Brazil, \EMAIL{vincent.guigues@fgv.br}} 
\AUTHOR{Anton J. Kleywegt}
\AFF{Industrial and Systems Engineering, Georgia Institute of Technology, Atlanta, USA,
\EMAIL{anton@isye.gatech.edu}}
\AUTHOR{Victor Hugo Nascimento}
\AFF{Systems Engineering and Computer Science, UFRJ, Rio de Janeiro, Brazil, \EMAIL{victorrn@cos.ufrj.br}}
} 

\ABSTRACT{%
We introduce new optimization models for the dispatch of ambulances in two types of situations.
In the first type of situation, an emergency call arrives, and a decision has to be made whether an ambulance should be dispatched for that call, and if so, which ambulance should be dispatched, or whether the request should be put in a queue of waiting requests.
In the second type of situation, an ambulance finishes its current task, and a decision has to be made whether the ambulance should be dispatched to a request waiting in queue, and if so, which request, or whether the ambulance should be dispatched to an ambulance staging location, and if so, which ambulance staging location.
These ambulance dispatch decisions affect not only the emergency call and ambulance under consideration, but also the ability of the ambulance fleet to respond to future emergencies.
There is uncertainty regarding the locations, arrival times, types, service needs, and service times of future emergencies, as well as uncertainty regarding travel times.
The types and service needs of emergencies are important, because most emergency medical services have crew members with different qualifications and skills and ambulances with different capabilities, and emergencies should be served by crews and ambulances appropriate for the type of emergency.
A number of methods have been proposed in the literature for making ambulance dispatch decisions, none of which takes into account the importance of dispatching a crew and ambulance with the right capabilities for the emergency at hand.
Also, it is well known that response time is important.
However, response time is not of equal importance for all emergencies, and sometimes there is a trade-off between dispatching the ambulance/crew with the shortest response time and dispatching the ambulance/crew with the best capabilities for the emergency.
The methods proposed in this paper are the first that take into account how well the capabilities of each ambulance/crew would serve the emergency at hand, and the trade-off between that and response time.
In addition, the methods take into account the impact of current dispatch decisions on the ability to respond to uncertain future emergencies.
We compare the dispatching policies proposed in this paper with policies previously proposed in the literature, using data of more than $2$~years of emergency calls for the Rio de Janeiro emergency medical service.
These tests show that the proposed policies result in both shorter response times as well as better matching of ambulances/crews with emergencies than the other policies.
}




\maketitle

%


\section{Introduction}
\label{sec:intro}

This paper proposes new optimization models for aiding ambulance dispatch decisions.
The paper also proposes solution methods for these problems, and compares the performance of the proposed methods with methods previously proposed in the literature, using data from the Rio de Janeiro emergency medical service.
Next we describe essential aspects of ambulance operations.

\subsection{Ambulance Fleet Operations}
\label{sec:ambulance operations}

We consider requests for emergency medical service that arrive at a call center.
A call center telecommunicator receives the call and obtains data from the caller, including the nature of the emergency and the location of the emergency.
The telecommunicator records the data, classifies the emergency, and decides whether to request that an ambulance be dispatched to the emergency.
If so, a decision is made which ambulance should be dispatched to the emergency, or whether the request should be placed in a queue of requests waiting for an ambulance to be dispatched.

There are multiple systems for classifying emergencies.
Most systems classify emergencies into 30--40 chief complaint types, based on the anatomical system involved (such as headache) or the cause of the emergency (such as burns/explosion), and the emergencies in each chief complaint type are further classified into subtypes.
The subtypes provide information regarding the importance of response time as well as the capabilities of the ambulance and crew that should treat such an emergency.
We call the combination of chief complaint type and subtype the {\em emergency type} or {\em call type}.

Usually, if an ambulance is dispatched to an emergency, it is an available ambulance, that is, an ambulance that is not on the way to an emergency, or busy at an emergency, or transporting patients from an emergency to a hospital.
That is, an available ambulance is an ambulance that is either waiting at an ambulance staging location for an assignment, or on its way to an ambulance staging location after attending to an emergency.
Some emergency medical services (EMSs) also consider dispatching an ambulance that is on the way to another emergency, that is, an ambulance that has already been dispatched to an (apparently less urgent) emergency is preempted and dispatched to the new (apparently more urgent) emergency.
In some systems the telecommunicator makes the ambulance selection decision, and in other systems a separate ambulance dispatcher makes the decision.

The amount of time that elapses from the moment the first call related to an emergency is received until the first ambulance personnel arrive at the patient(s) is called {\em the response time}.
(In practice, distinction is made between different response times, for example, the elapsed time can start when the first call related to an emergency is received or when the telecommunicator requests that an ambulance be dispatched or when the ambulance personnel receive the dispatch instructions, and the elapsed time can end when the first ambulance arrives at the location of the emergency or when the first ambulance personnel arrive at the patient(s).
In this paper we consider the first of these times as the response time.)
It has been found that the probability of survival to hospital admission and the probability of survival to hospital discharge often decrease as response time increases, but the dependence varies greatly by emergency type and by treatment, and sometimes the dependence is hard to detect statistically \citep{cret:79,blac:91,lars:93,vale:97,stie:99,vale:00,waal:01,pell:01,dema:03,pons:05,blac:09,blan:12,weis:13}.
In addition, response times of EMSs are relatively easy to measure, and most EMSs record and report (for example, to the National Emergency Medical Services Information System NEMSIS) response time data.
Therefore many EMSs as well as academic papers put great emphasis on the response time performance metric.
As summary statistic of response time data, EMSs often measure specific quantiles of the response time empirical distribution.
For example, the $0.8$ and $0.9$ empirical quantiles are measured (often pooling all emergency data, that is, without distinction between different types of emergencies), and for performance to be called acceptable, the considered quantiles have to be less than specified threshold values \citep{hend:04,rest:09,maso:13}.

It has been pointed out that response time is not the only factor under control of EMSs that affects the survival probabilities of patients, and also that the impact of emergency medical care and response time greatly depends on the emergency type.
Also, the capabilities of the ambulance and personnel can affect the survival probability of the patient, depending on the emergency type.
Although many academic papers consider all ambulances to be the same, typically ambulances are not the same.
For example, some EMSs make a distinction between basic life support (BLS) ambulances and advanced life support (ALS) ambulances \citep{schi:79}.
Some EMS systems such as the one in Rio de Janeiro (BLS, Intermediate Life Support, and ALS), and the Grady system in Atlanta (BLS, ALS, and Stroke Unit), have more than two types of ambulances.
Regarding personnel, many EMS systems have Emergency Medical Technician (EMT) and Paramedic personnel.
Some EMSs make a distinction between first responders staffed by EMTs and second responders staffed by Paramedics.
Some EMS systems also have Advanced Emergency Medical Technician (AEMT) personnel, and some have crews that include nurses and physicians with various specialties.
At a sufficient level of detail, every ambulance is unique, because the crew members of different ambulances have different qualifications and experience.
\cite{pers:03} compared a targeted EMS system that dispatches BLS or ALS ambulances according to the emergency type and an EMS system that provides uniform service with all ALS ambulances.
The study focused on witnessed ventricular fibrillation cardiac arrest emergencies, and compared system performance in terms of four outcomes: return of spontaneous circulation, survival to hospital admission, survival to hospital discharge, and survival to 1 year.
All performance measures were better for the targeted system that dispatched BLS or ALS ambulances according to the emergency type.

A potential concern with using the emergency type in dispatch decisions is that emergencies are sometimes misclassified by telecommunicators due to limited or incorrect data.
Based on an analysis of medical emergency data in the UK under two systems for classifying emergencies, \cite{nich:99} found that the benefits of dispatch based on emergency type outweighed the risk of misclassification.
Therefore, it makes sense to take the emergency type into account when deciding whether to immediately dispatch an ambulance to the emergency, and if so, which ambulance to dispatch, as opposed to considering only the response time irrespective of the emergency type.
For the reasons discussed above, our models take the type of each emergency as well as the capabilities of each individual ambulance into account when evaluating ambulance dispatch decisions.

After arriving at the location of the emergency, the ambulance crew treat the patients, and decide whether to transport patients to a hospital.
In some cases, it is decided not to transport any patients to a hospital, and thus the ambulance becomes available after treatment at the emergency site has been completed.
In other cases, the ambulance transports patients to a chosen hospital, and becomes available thereafter.
Sometimes ambulances are used for scheduled transport of patients, for example, to transfer patients from one hospital to another.
Such transport tasks are regarded as a specific ``emergency'' type.
After patients have been transported in an ambulance, the ambulance has to be cleaned.
The cleaning is regarded as part of the ambulance's work related to an emergency.
After an ambulance has completed its work related to an emergency, a decision has to be made where to dispatch the newly available ambulance to.
If any requests are waiting in queue, the ambulance can be dispatched to a chosen waiting emergency.
The ambulance can also be used to preempt another ambulance that is on its way to the location of an emergency, for example, if the newly available ambulance is closer to the emergency than the previously dispatched ambulance.
Otherwise the ambulance can be sent toward an ambulance staging location, and it is regarded as available during such a journey.
There are various types of ambulance staging locations.
Most EMSs have one or more facilities that make provision for ambulance maintenance, office space, and training.
Such facilities can be used for ambulance staging.
Often, a larger number of smaller facilities are also used for ambulance staging.
These smaller facilities may make provision for ambulance personnel to relax during less busy times.
During busy times, ambulances may use public facilities such as parking lots for ambulance staging.
For the purpose of this paper, all ambulance staging locations will be called stations.

\subsection{Contributions}
\label{sec:contributions}

The contributions of this paper are as follows.

\noindent
\paragraph{\textbf{(A) Modeling the operations of an ambulance fleet, incorporating important features omitted in previous publications.}}
We propose optimization models for ambulance dispatch decisions that include the following features:
\begin{enumerate}
\item
The models take into account the capabilities of each ambulance and crew.
The capabilities vary among ambulances and especially among crews.
For example, paramedics are allowed to administer various medications that EMTs are not allowed to administer, some crews have received special training for specific emergencies such as the handling of stroke victims, and some crews are more skilled in coping with dangerous situations such as rioting.
As pointed out in Section~\ref{sec:literature}, none of the existing methods for ambulance dispatch decisions take into account differences in the capabilities of different ambulances and crews.
\item
The models take into account the type of each emergency.
As mentioned in Section~\ref{sec:ambulance operations}, the emergency type provides information regarding the anatomical system involved, the cause of the emergency, the importance of response time, as well as the capabilities of the ambulance and crew that should treat such an emergency.
Thus emergency type contains more information than just a priority level, and this information is important for determining the type of ambulance and crew needed, the marginal value of response time, and the set of appropriate hospitals for the emergency.
As pointed out in Section~\ref{sec:literature}, existing methods for ambulance dispatch decisions either consider only one emergency type, that is, the effects of emergency type mentioned above are ignored, or a small number ($2$ or $3$) of priority levels.
\item
The models allow the hospital for each emergency to be chosen, taking into account the emergency type and the location of the emergency relative to hospitals.
Existing methods for ambulance dispatch decisions ignore choices among multiple hospitals; in fact, most existing models either have no hospital entity in the model, or have a single ``hospital'' entity irrespective of the type and location of the emergency.
\item
The models make provision for a queue of waiting emergencies.
This is both necessary if all ambulances are busy, and desirable if the emergency is not urgent and few ambulances are available.
As pointed out in Section~\ref{sec:literature}, many existing models ignore the possibility of a queue of waiting emergencies.
Instead, those models assume that if all ambulances are busy then some unlimited outside service will take care of the emergency.
\item
The models make provision for both ambulance selection decisions as well as nontrivial ambulance reassignment decisions.
When an ambulance becomes available, it can be assigned to an emergency in queue, and it can also be sent to a chosen station.
Since many existing models ignore the possibility of a queue of waiting emergencies, such models cannot accommodate decisions to assign available ambulances to emergencies in queue.
Also, many existing models assume that when an ambulance becomes available, it will go to its home base, even if sending it to a different station would improve coverage greatly.
\item
The models allow ambulances on their way to a station to be dispatched to an emergency.
Because of the triangle inequality, that results in smaller response times than waiting until the ambulance reaches the station before dispatching the ambulance to the emergency.
Also, with modern communication technology dispatchers are tracking ambulances and are in constant contact with ambulances, and such en-route dispatching is easy to execute.
As pointed out in Section~\ref{sec:literature}, most existing models allow only ambulances at stations to be dispatched.
\item
Dispatch decisions have consequences not only for the emergency and the ambulance under consideration, but also for future emergencies and for other ambulances that have to take care of future emergencies.
It is challenging to take these future consequences of dispatch decisions into account, because future consequences are a complicated function of current decisions, and because the future consequences are uncertain.
The models take into account uncertain future consequences.
As pointed out in Section~\ref{sec:literature}, existing models either ignore future consequences, or use an indirect quantity such as a ``coverage'' metric or a ``preparedness'' metric to make some provision for future consequences.
\end{enumerate}

\noindent
\paragraph{\textbf{(B) Practical and competitive methods for making ambulance dispatch decisions.}}
This paper proposes two optimization-based methods for each of the two types of decisions mentioned above:
\begin{enumerate}
\item
When a request arrives, the decision whether to dispatch an ambulance to the emergency immediately, and if so, which ambulance to dispatch to the emergency, or whether to add the request to a queue of waiting requests.
This decision is called the {\em ambulance selection decision}.
\item
When an ambulance becomes available (after completing its service at the location of an emergency or at a hospital), the decision whether to dispatch the ambulance to an emergency waiting in queue, and if so, to which emergency to dispatch the ambulance, or whether to send the ambulance to a station, and if so, to which station to send the ambulance.
This decision is called the {\em ambulance reassignment decision}.
\end{enumerate}

We do not include the two types of preemption alternatives mentioned above, that is, (1)~when a request arrives, the alternative to preempt an ambulance already on its way to an emergency and to dispatch it to the newly arrived emergency, and (2)~when an ambulance becomes available, the alternative to preempt an ambulance already on its way to an emergency and to dispatch the newly available ambulance to the emergency (called ``diversions'' by \cite*{maso:13}).
\cite{lim:11} simulated and compared four dispatch policies, obtained by switching each of the two preemption/diversion alternatives mentioned above on and off (the first alternative was called ``Reroute Enabled Dispatch'' and the second alternative was called ``Free Ambulance Exploitation Dispatch'').
The results suggested that these preemption/diversion alternatives add very little benefit.
In addition, preempting ambulances have practical disadvantages, and thus it is questionable whether use of these alternatives is a wise decision.

For each type of ambulance dispatch decision, we propose two optimization models.
Each of these models incorporate all the problem features discussed in (A) above, including multiple emergency types, ambulance and crew capabilities, hospital selection, queues of waiting emergencies, en-route dispatch, and uncertain future consequences of dispatch decisions.
We show that even for large EMS systems such as the one in Rio de Janeiro, these optimization problems can be solved in at most a few seconds --- fast enough for use as a practical dispatch decision support tool.
In addition we compare the performance of the two methods with methods that have been proposed in the literature for ambulance dispatch decisions, in terms of both the traditional response time performance metric as well as a weighted response time performance metric that takes into account the emergency type and the capabilities of the ambulance and crew dispatched to the emergency.
We used more than 2 years of emergency call data of the Rio de Janeiro EMS, and based on these data showed that the two proposed methods consistently outperform the other methods.

The rest of this paper is organized as follows.
Section~\ref{sec:literature} reviews the literature on operations research work for ambulance applications, with emphasis on work to support ambulance dispatch decisions.
In Section~\ref{sec:optmodel}, we specify the optimization models.
In Section~\ref{sec:sol}, we explain how these problems are solved.
The numerical results for the Rio de Janeiro EMS are presented in Section~\ref{sec:num}.

\section{Related Literature}
\label{sec:literature}

Most of the literature related to ambulance operations address the location and relocation of emergency facilities, including stations and ambulances.
For surveys of this literature, see \cite{brot:03,gree:04,lizh:11,arin:17}.
Compared with the work on static location problems for stations and/or ambulances, relatively little work has been done to optimize ambulance operations.
In ambulance operations, two types of dispatch decisions mentioned in Section~\ref{sec:ambulance operations} are important: (1)~the ambulance selection decision, and (2)~the ambulance reassignment decision.
The solutions of the static location problems are sometimes used for ambulance reassignment decisions by making a static assignment of a home base to each ambulance, and whenever an ambulance becomes available and is not dispatched to a request waiting in queue, the ambulance is sent to its home base \citep{gold:90b,hend:04,rest:09,band:12,knig:12,maso:13,mayo:13,band:14}.
A dynamic extension of this approach for ambulance reassignment is to choose which station to send an ambulance to each time when an ambulance becomes available and is not dispatched to a request waiting in queue.
This dynamic reassignment approach is followed in \cite{schm:12} and in this paper.

For the ambulance selection decision, the closest available ambulance rule is simple and popular \citep{berl:74,gold:90a,hend:99,hend:04,maxw:09,maxw:10,maxw:13,alan:13}.
Some EMSs partition the service region into response areas or districts, and apply some fixed preference policy to make ambulance selection decisions.
For example, for each district, a preference list of districts is chosen in advance.
Typically, for each district, that district appears first in its preference list.
Then, when an emergency call located in a district arrives, the first district in the emergency district's priority list with an available ambulance is determined, and the ambulance in that district closest to the emergency location is dispatched.
\cite{cart:72} conducted a detailed study for a setting with 2 ambulances, and characterized the optimal response area for each ambulance.
\cite{swov:73a} compared two dispatch rules, the closest available ambulance rule, and a service district rule that works as follows:
Each ambulance is assigned to a service district, with one ambulance assigned to each district.
When an emergency call arrives in a district, if the ambulance assigned to that district is available, then it is dispatched to the call, even if it is temporarily outside its district; otherwise, if the ambulance assigned to that district is not available, then the closest available ambulance is dispatched to the call.
\cite{knig:12} proposed a dispatch rule based on a static preference matrix $\rho$.
The service region is partitioned into ``demand nodes'', and each available ambulance is assigned to a station.
Then $\rho_{i,j}$ denotes the $j$th most preferred station to use for an emergency at demand node~$i$.
An ambulance is dispatched from station $\rho_{i,j}$ to an emergency at demand node~$i$ if and only if there is no available ambulance at stations $\rho_{i,1},\ldots,\rho_{i,j-1}$ and there is at least one available ambulance at station $\rho_{i,j}$.

\cite{ande:07} proposed heuristics for ambulance dispatch and relocation for a system with three priority levels, based on a measure of ``preparedness'' for each zone.
For priority~1 calls, they dispatch the closest available ambulance.
For priority~2 and~3 calls, they dispatch the ambulance with expected travel time less than a specified threshold that will result in the least decrease in the minimum preparedness measure over all zones.
They did not specify how they select a request waiting in queue for the ambulance reassignment decision.
\cite{lees:11} considered the same preparedness measure as \cite{ande:07}, and showed that it resulted in worse performance than the closest-available-ambulance rule.
Then two modifications of the preparedness-based dispatching rule were proposed.
The first modification dispatches the available ambulance that maximizes the minimum preparedness measure over all zones divided by the travel time from the ambulance to the emergency location.
The second modification replaces the minimum preparedness measure over all zones in the calculations with other aggregates of the preparedness measures of different zones.
If an ambulance becomes available and there are requests waiting in queue, then the newly available ambulance is dispatched to the request in queue closest to the ambulance.
\cite{lees:12} proposed a rule for the ambulance reassignment decision that takes into account both the distances or times between the newly available ambulance and requests waiting in queue, as well as a centrality measure of each request waiting in queue.
It was not specified how to select a station to which to send the newly available ambulance if there are no requests waiting in queue, or how the ambulance selection decision is made.

\cite{schm:12} proposed a dynamic programming formulation for the ambulance selection decision and the ambulance reassignment decision, with the restriction that if an ambulance becomes available and there are requests waiting in queue, then the newly available ambulance is dispatched to the next request in first-come-first-served order.
An approximate dynamic programming method was used that restricted ambulance selection decisions to the closest-available-ambulance rule, and that attempted to optimize the ambulance reassignment decision.
It was shown that after sufficient training these solutions outperform policies that combine the closest-available-ambulance dispatching rule with one of the following three ambulance reassignment rules if there are no requests waiting in queue: (1)~send the newly available ambulance to the closest station, (2)~send the newly available ambulance to its home base, or (3)~send the newly available ambulance to a random station.

\cite{band:12} formulated an ambulance dispatching problem with two priority levels and exponentially distributed service times as a continuous-time Markov Decision Process (MDP), and showed (for a sufficiently small number of ambulances to enable solving the MDP) that the closest-available-ambulance dispatching policy is suboptimal.
It was assumed that if an emergency call arrives and all ambulances are busy, then the emergency is handled by an ``outside'' service, for example provided by the fire department, that is, there is no queue of waiting requests in the model.
It was also assumed that an ambulance returns to its home base when it becomes available.
Similarly, \cite{band:14} proposed an ambulance dispatching heuristic that takes emergency priorities into account, and used simulation to compare the performance of the heuristic and the closest-available-ambulance dispatching rule for a setting with two priority levels.
The same assumptions as in \cite{band:12} are made, including the assumption that an ambulance returns to its home base when it becomes available.

\cite{mayo:13} proposed a heuristic to partition the service region into districts, locate stations in each district, and assign each ambulance to a station.
They used simulation to compare the performance of four static dispatch policies for a setting with two priority levels; two types of policies specifying ambulance selection decisions if there is an ambulance available in the same district as the emergency, combined with two types of policies specifying ambulance selection decisions if there is no ambulance available in the same district.
If there is an ambulance available in the same district as the emergency, then the first type of policy dispatches the closest available ambulance within the same district, and the second type of policy applies a heuristic of \cite{band:12} to each district.
If there is no ambulance available in the same district as the emergency, then the first type of policy assumes that an outside emergency response is automatically dispatched within the same district, and the second type of policy dispatches an ambulance from another district using a preference list of ambulances.
For all policies, it is assumed that if an emergency call arrives and all ambulances are busy, then the emergency is handled by an outside service, that is, no queue of waiting requests is considered.
Also, it was assumed that an ambulance returns to its home base when it becomes available.

\cite{lisay:16} used simulation to compare two dispatch policies, a policy that dispatches the closest available ambulance to all calls, and a policy that dispatches the closest available ambulance to priority~1 calls, and the ambulance within a specified response time radius which has the least utilization to priority~2 and~3 calls.
It was assumed that if an emergency call arrives and all ambulances are busy, then the call is lost.
It was also assumed that an ambulance returns to its home base when it becomes available.

\cite{jagt:17a} compared two dispatch policies with the closest-available-ambulance policy.
In all policies, ambulance reassignment decisions are made according to the first-come first-served rule: if an ambulance becomes available and there are calls waiting in queue, then the ambulance is dispatched to the waiting call that first entered the queue.
If an ambulance becomes available and there is no call waiting in queue, then the ambulance returns to its home base.
One dispatch policy is based on a Markov decision process with a simplified state, that represents the location of the currently considered emergency, and the set of available ambulances, assuming that each available ambulance is at its home base.
The second dispatch policy, called DMEXCLP, works as follows: for each available ambulance that can reach the emergency location within the threshold time, the objective value of the maximum expected covering location problem (MEXCLP), originally proposed by \cite{dask:83}, with the available ambulances excluding that ambulance, is computed.
Then the available ambulance that can reach the emergency location within the threshold time with the largest objective value of the MEXCLP without that ambulance is dispatched.
If no available ambulance can reach the emergency location within the threshold time, then the available ambulance, irrespective of travel time to the emergency location, with the largest objective value of the MEXCLP without that ambulance is dispatched.
Simulation results showed that the heuristic has a much lower fraction of late arrivals than the closest-available-ambulance policy, but that the heuristic also has a much greater mean response time than the closest-available-ambulance policy.
\cite{jagt:17b} assumed that all ambulances are dispatched from their home bases, that there is always an ambulance available at its home base when an emergency call arrives, and that an ambulance must always be dispatched immediately to an emergency.
They also assumed that after service each ambulance returns to its home base, and that the time that elapses from the moment the ambulance arrives at the emergency location until the ambulance is back at its home base is deterministic and is the same for all emergencies and all ambulances, that is, the elapsed time does not depend on the emergency location or the ambulance home base location.
They showed with an example that for any online policy the ratio of fraction of late arrivals to the offline minimum fraction of late arrivals can be arbitrarily large, and thus no online policy can have a finite competitive ratio.
They also used simulation to compare the expected performance ratios of the closest-available-ambulance policy and of DMEXCLP, and showed that the expected performance ratio of DMEXCLP is better than that of the closest-available-ambulance policy.

\section{Optimization Models}
\label{sec:optmodel}

We consider the following two types of ambulance dispatch problems:
\begin{enumerate}
\item
The ambulance selection problem:
When a request arrives, the problem to decide whether an ambulance should be dispatched to the emergency, and if so, which ambulance to dispatch to the emergency, or whether to add the request to a queue of waiting requests.
\item
The ambulance reassignment problem:
When an ambulance becomes available, the problem to decide whether to dispatch the ambulance to an emergency waiting in queue, and if so, to which emergency to dispatch the ambulance, or whether to send the ambulance to a station, and if so, to which station to send the ambulance.
\end{enumerate}
For each of these two types of ambulance dispatch problems, we propose two optimization models, one ``itinerary-based'' optimization model and one ``arc-based'' optimization model.
All the optimization models are intended to be solved in rolling horizon fashion, each time that an ambulance selection decision or an ambulance reassignment decision has to be made.
All models ``look ahead'' until the end of a chosen time horizon (for example, a few hours or until the end of the day) to approximate the impact of current decisions on the objective values in the future.
All problems minimize a combination of the cost of the immediate decision and the expected future costs affected by the immediate decision over the planning horizon.
Each time that an ambulance dispatch decision has to be made, a two-stage stochastic optimization problem is solved.
The first-stage decisions are either ambulance selection decisions or ambulance reassignment decisions, as appropriate for the situation at hand.
The second-stage decisions are sequences of ambulance selection decisions and ambulance reassignment decisions over the planning horizon.
For each first-stage decision, the setting for the decision is known.
For example, if an ambulance selection decision has to be made then the type and location of the newly arrived emergency is known, or if an ambulance reassignment decision has to be made then the location of the newly available ambulance is known, and in both cases the current state of the system is known.
In contrast, the settings for second-stage decisions are not known yet, and the uncertainty is represented with a set of second-stage scenarios.
The number of first-stage alternatives is relatively small.
If an ambulance selection decision has to be made, then the first-stage alternatives correspond to the available compatible ambulances combined with the candidate hospitals, as well as the alternative to add the emergency to the queue.
If an ambulance reassignment decision has to be made, then the first-stage alternatives correspond to the emergencies in queue as well as the stations.
Therefore, each two-stage problem can be solved by computing the expected second-stage cost for each feasible first-stage alternative, and then choosing the first-stage alternative with the best combination of first-stage and second-stage cost.
Therefore, to solve the two-stage problems fast enough, we need to be able to quickly compute the expected second-stage cost for every feasible first stage decision.
The itinerary-based and the arc-based models differ in the way that solutions for the second-stage problems are represented.
The itinerary-based model represents a second-stage solution with a single itinerary for each ambulance consisting of a sequence of dispatches/tasks, whereas the arc-based model represents a second-stage solution by stringing together a sequence of separate dispatch decisions for each ambulance.

Next we present the features that are common to all the models.
Additional features that are specific to some models are added later.
Ambulances can be dispatched from stations, from emergency locations (if the ambulance is not needed to transport patients to a hospital), from hospitals, and from intermediate locations while traveling towards a station.
The models considered here do not allow ambulances to be dispatched while busy with service --- while traveling towards an emergency location, or while providing on-site emergency care, or while traveling with patient(s) towards a hospital.
That is, we do not model preemption or ``forward'' dispatching of ambulances.
An ambulance that is not busy with service must either be at a station, or traveling toward a station.
If ambulances are allowed to wait at a hospital for a dispatch, then the hospital is also a station in the model.

Let $[0,T]$ denote the time horizon currently considered, with $t=0$ denoting the current time.
Thus the dispatch decision on-hand has to be made at time~$0$.
We assume that it does not happen simultaneously that an emergency call arrives and an ambulance becomes available.
Therefore, at (each) time $t=0$, either an ambulance selection problem or an ambulance reassignment problem is considered.

The feasible decisions at time~$0$ are the same for the itinerary-based ambulance selection problem and the arc-based ambulance selection problem, and similarly the feasible decisions at time~$0$ are the same for the itinerary-based ambulance reassignment problem and the arc-based ambulance reassignment problem.
Space is discretized for all the models.
Let $\mathcal{L}$ denote the set of discrete locations, used for representing emergency call locations as well as ambulance locations.
Each emergency call is characterized by its arrival time, its location, and its type.
Let $\mathcal{C}$ denote the set of emergency types,
let $\mathcal{A}$ denote the set of ambulances and their crews,
let $\mathcal{B}$ denote the set of ambulance stations,
and let $\mathcal{H}$ denote the set of hospitals.
An emergency of type $c \in \mathcal{C}$ can be served by a subset $\mathcal{A}(c) \subset \mathcal{A}$ of ambulances, and an emergency of type $c \in \mathcal{C}$ at location $\ell \in \mathcal{L}$ can be sent to a subset $\mathcal{H}(c,\ell) \subset \mathcal{H}$ of hospitals.
Many emergency types do not involve transport to a hospital.
To simplify notation, it is assumed that the emergency type $c \in \mathcal{C}$ contains information whether transport to a hospital is involved or not, and if emergency type~$c$ does not involve transport to a hospital, then $\mathcal{H}(c,\ell) = \ell$.
Thus the set $\mathcal{H}$ of ``hospitals'' contains more than actual hospitals.
Let $\mathcal{C}(a) \defi \{c \in \mathcal{C} \; : \; a \in \mathcal{A}(c)\}$ denote the set of emergency types that can be served by ambulance~$a$.

The initial conditions at time~$t=0$ are specified as follows:
If a call arrives at time~$0$, then the location and the type of the emergency that has just arrived are denoted by $\ell_{0}$ and $c_{0}$ respectively.
If an ambulance completes service (at a ``hospital'') at time~$0$, then the ambulance and the hospital are denoted by $a_{0}$ and $h_{0}$ respectively.
Let $C_{0}(c,\ell)$ denote the number of emergencies of type~$c$ at location~$\ell$ waiting in queue at time~$0$ for an ambulance to be dispatched to serve the emergency, including any call of type~$c$ at location $\ell$ that has just arrived at $t=0$.
Let $A_{0}(a,b) \in \{0,1\}$ indicate whether ambulance~$a$ is available at~$b$ for dispatch just before the dispatch of ambulances at time~$0$.

First we present deterministic formulations of the problems, as though the arrival times, locations, types of the emergency calls, travel times, and emergency service times over the planning horizon are known.
In Section~\ref{sec:smodel} we describe the extension of this model to incorporate random arrival times, locations, types of emergency calls, travel times, and emergency service times.

\subsection{Deterministic Arc-based Models}
\label{sec:detarcmodel}

\subsubsection{Problem input parameters}
\label{sec:parameters}

Time is discretized for the model, with $t=0$ denoting the current time, and $t = 1,\ldots,T$ denoting the time steps until the end $T$ of the time horizon.

\textbf{Keeping track of ambulance locations.}
As mentioned above, ambulances can be dispatched while traveling towards a station.
To model this, and in general to keep track of ambulance locations both while stationary and while moving, it is useful to determine where ambulances can be while traveling to specific destinations.
First,
let $L(t,a,h,b) \in \mathcal{L}$ denote the location of ambulance~$a$ at time~$t+1$ if the ambulance starts from hospital~$h$ at time~$t$ to travel towards station~$b$.
In addition,
let $L(t,a,\ell_{1},b) \in \mathcal{L}$ denote the location of ambulance~$a$ at time~$t+1$ if the ambulance is at location $\ell_{1}$ at time~$t$ and continues to travel towards station~$b$.
Next, for each time $t \in \{1,\ldots,T\}$, let
\[
\mathcal{L}_{1}(t,a,b) \ \ \defi \ \ \bigcup_{h \in \mathcal{H}} \{L(t-1,a,h,b)\} \setminus \{b\}
\]
denote the set of all intermediate locations that can be reached by an ambulance at time~$t$ if the ambulance starts traveling at time~$t-1$ from some hospital towards~$b$.
By induction, for each $\tau \ge 2$, let
\[
\mathcal{L}_{\tau}(t,a,b) \ \ \defi \ \ \bigcup_{\ell_{1} \in \mathcal{L}_{\tau-1}(t-1,a,b)} \{L(t-1,a,\ell_{1},b)\} \setminus \{b\}
\]
denote the set of all intermediate locations that can be reached by an ambulance at time~$t$ if the ambulance starts traveling at time~$t-\tau$ from some hospital towards~$b$.
Then,
let
\[
\mathcal{L}(t,a,b) \ \ \defi \ \ \bigcup_{\tau \ge 1} \mathcal{L}_{\tau}(t,a,b)
\]
denote the set of all intermediate locations that can be reached by ambulance~$a$ at time~$t$ if the ambulance travels from some hospital towards station~$b$.

Let $A_{0}(a,\ell_{1},b) \in \{0,1\}$ indicate whether ambulance~$a$ is at location~$\ell_{1}$ traveling towards station~$b$ available for dispatch just before the dispatch of ambulances at time~$0$.
In addition,
let $A_{0}(c,a,\ell_{1},\ell,h) \in \{0,1\}$ indicate whether ambulance~$a$ is at location~$\ell_{1}$ at time~$0$ traveling to an emergency type~$c$ at location~$\ell$ and from there to hospital~$h$.
Also, let $A_{0}(c,a,\ell_{1},h) \in \{0,1\}$ indicate whether ambulance~$a$ is at location~$\ell_{1}$ at time~$0$ traveling to hospital~$h$ with emergency type~$c$ patient(s) after on-site emergency care has already been provided.

\textbf{Future emergency calls, service times, and travel times.}
For each time $t \in \{1,\ldots,T\}$, emergency type $c \in \mathcal{C}$, and emergency location $\ell \in \mathcal{L}$, let $\lambda(t,c,\ell)$ denote the number of calls of type~$c$ at location~$\ell$ arriving in period~$t$.
Let $\tau(t,c,a,\ell_{1},\ell,h)$ denote the time for ambulance~$a$ to travel from $\ell_{1}$ at time~$t$ to $\ell$, provide on-site emergency care for emergency type~$c$ at $\ell$, travel with patient(s) from $\ell$ to hospital~$h$, and deliver the patient(s) at~$h$.
Also,
let $\tau_{0}(c,a,\ell_{1},h)$ denote the time for ambulance~$a$ to travel with patient(s) from $\ell_{1}$ at time~$0$ after on-site emergency care has already been provided, to hospital~$h$, and deliver the patient(s) of type~$c$ at $h$.

\subsubsection{Decision variables for the ambulance selection problem}

The following first-stage (for $t=0$) decision variables are used for the ambulance selection problem:
\begin{itemize}
\item
$x_{0}(c_{0},a,b,\ell_{0},h) \in \{0,1\}$ indicates whether ambulance~$a \in \mathcal{A}(c_{0})$ is dispatched at time~$0$ from station $b \in \mathcal{B}$ to serve the emergency of type $c_{0}$ that arrived at location $\ell_{0}$ at time~$0$, and transport patient(s) to hospital $h \in \mathcal{H}(c_{0},\ell_{0})$ (such a decision variable is used only if $A_{0}(a,b) = 1$);
\item
$x_{0}(c_{0},a,\ell_{1},b,\ell_{0},h) \in \{0,1\}$ indicates whether ambulance~$a \in \mathcal{A}(c_{0})$ at location $\ell_{1} \in \mathcal{L}(0,a,b)$ at time~$0$ traveling toward station $b \in \mathcal{B}$ is dispatched to serve the emergency of type~$c_{0}$ that arrived at location $\ell_{0}$ at time~$0$, and transport patient(s) to hospital $h \in \mathcal{H}(c_{0},\ell_{0})$ (such a decision variable is used only if $A_{0}(a,\ell_{1},b) = 1$).
\end{itemize}

The following second-stage (for $t=1,\ldots,T$) decision variables are used for the ambulance selection problem:
\begin{itemize}
\item
$x_{t}(c,a,b,\ell,h) \in \{0,1\}$ indicates whether ambulance~$a \in \mathcal{A}(c)$ is dispatched at time~$t$ from station $b \in \mathcal{B}$ to serve an emergency of type $c \in \mathcal{C}$ (including both emergencies in queue as well as any calls that arrived at time~$t$) at location $\ell \in \mathcal{L}$, and transport patient(s) to hospital $h \in \mathcal{H}(c,\ell)$;
\item
$x_{t}(c,a,\ell_{1},b,\ell,h) \in \{0,1\}$ indicates whether ambulance~$a \in \mathcal{A}(c)$ at location $\ell_{1} \in \mathcal{L}(t,a,b)$ at time~$t$ traveling toward station $b \in \mathcal{B}$ is dispatched to serve an emergency of type $c \in \mathcal{C}$ that arrived at location $\ell \in \mathcal{L}$ and transport patient(s) to hospital $h \in \mathcal{H}(c,\ell)$;
\item
$x_{t}(c,a,h',\ell,h) \in \{0,1\}$ indicates whether ambulance~$a \in \mathcal{A}(c)$ is dispatched at time~$t$ from hospital $h' \in \mathcal{H}$ to serve an emergency of type $c \in \mathcal{C}$ at location $\ell \in \mathcal{L}$, and transport patient(s) to hospital $h \in \mathcal{H}(c,\ell)$;
\item
$y_{t}(a,h,b) \in \{0,1\}$ indicates whether ambulance~$a \in \mathcal{A}$ is instructed at time~$t$ to move from hospital $h \in \mathcal{H}$ towards station $b \in \mathcal{B}$;
\item
$C_{t}(c,\ell) = $ the number of emergencies of type~$c \in \mathcal{C}$ waiting in queue at location $\ell \in \mathcal{L}$ at the beginning of period~$t$;
\item
$A_{t}(a,b) \in \{0,1\}$ indicates whether ambulance~$a \in \mathcal{A}$ is at station $b \in \mathcal{B}$ at the beginning of period~$t$;
\item
$A_{t}(a,\ell_{1},b) \in \{0,1\}$ indicates whether ambulance~$a \in \mathcal{A}$ is at location $\ell_{1} \in \mathcal{L}(t,a,b)$ moving towards station $b \in \mathcal{B}$ at the beginning of period~$t$.
\end{itemize}
Recall that the decision variables above are written for a deterministic problem; similar variables are used for a stochastic problem.

\subsubsection{Decision variables for the ambulance reassignment problem}

The following first-stage decision variables are used for the ambulance reassignment problem:
\begin{itemize}
\item
$x_{0}(c,a_{0},h_{0},\ell,h') \in \{0,1\}$ indicates whether ambulance~$a_{0}$ is dispatched at time~$0$ from hospital~$h_{0}$ to serve an emergency of type $c \in \mathcal{C}(a_{0})$ at location $\ell \in \mathcal{L}$, and transport patient(s) to hospital $h' \in \mathcal{H}(c,\ell)$ (such a decision variable is used only if $C_{0}(c,\ell) > 0$);
\item
$y_{0}(a_{0},h_{0},b) \in \{0,1\}$ indicates whether ambulance~$a_{0}$ is instructed at time~$0$ to move from hospital $h_{0}$ towards station $b \in \mathcal{B}$.
\end{itemize}

The second-stage (for $t=1,\ldots,T$) decision variables for the ambulance reassignment problem are the same as the decision variables for the ambulance selection problem.

\subsubsection{Constraints for the ambulance selection problem}

The following six sets of constraints apply to the first-stage variables for the ambulance selection problem only: \\

\paragraph{\textbf{(S1) Flow balance equations at the stations:}}
For each $a \in \mathcal{A}$, $b \in \mathcal{B}$,
\begin{align}
A_{1}(a,b) \ \ = \ \ & A_{0}(a,b) - \displaystyle{\mathds{1}_{\{a \in \mathcal{A}(c_{0})\}} \sum_{h \in \mathcal{H}(c_{0},\ell_{0})} x_{0}(c_{0},a,b,\ell_{0},h)} \nonumber \\
& {} + \displaystyle{\sum_{\{\ell_{1} \in \mathcal{L}(0,a,b) \, : \, L(0,a,\ell_{1},b) = b\}} \left[A_{0}(a,\ell_{1},b) - \mathds{1}_{\{a \in \mathcal{A}(c_{0})\}} \sum_{h \in \mathcal{H}(c_{0},\ell_{0})} x_{0}(c_{0},a,\ell_{1},b,\ell_{0},h) \right]}.
\label{eqn:flowbase0a}
\end{align}
The left side in~\eqref{eqn:flowbase0a} indicates whether ambulance~$a$ is at station~$b$ at the beginning of period $t=1$.
The first term $A_{0}(a,b)$ on the right indicates whether ambulance~$a$ is at station~$b$ at the beginning of period $t=0$.
The second term on the right indicates whether ambulance~$a$ leaves station~$b$ in the first stage to attend the call that has just arrived.
The third term on the right indicates whether ambulance~$a$ is scheduled to arrive at station~$b$ during the first stage.
The fourth term on the right indicates whether ambulance~$a$ was scheduled to arrive at station~$b$ during the first stage based on its status at $t=0$ but then is dispatched while en-route to station~$b$ to attend the call that has just arrived.
The remaining flow constraints follow the same logic and therefore are given without detailed explanation.

\paragraph{\textbf{(S2) Flow balance equations at the locations between hospitals and stations:}}
For each $a \in \mathcal{A}$, $b \in \mathcal{B}$, $\ell_{1} \in \mathcal{L}(0,a,b)$,
\begin{equation}
A_{1}(a,\ell_{1},b) \ \ = \ \ \displaystyle{\sum_{\{\ell'_{1} \in \mathcal{L}(0,a,b) \; : \; L(0,a,\ell'_{1},b) = \ell_{1}\}} \left[\displaystyle{A_{0}(a,\ell'_{1},b) - \mathds{1}_{\{a \in \mathcal{A}(c_{0})\}} \sum_{h \in H(c_{0},\ell_{0})} x_{0}(c_{0},a,\ell'_{1},b,\ell_{0},h)}\right]}.
\label{eqn:flowlocation0a}
\end{equation}

\paragraph{\textbf{(S3) Flow balance equations for the queues:}}
For each $c \in \mathcal{C}$, $\ell \in \mathcal{L}$,
\begin{align}
C_{1}(c,\ell) \ \ = \ \ & C_{0}(c,\ell) - \mathds{1}_{\{c = c_{0}, \, \ell = \ell_{0}\}} \left[\displaystyle{\sum_{a \in \mathcal{A}(c)} \sum_{b \in \mathcal{B}} \sum_{h \in \mathcal{H}(c,\ell)} x_{0}(c,a,b,\ell,h)}\right. \nonumber \\
& \left.\hspace{40mm} {} + \displaystyle{\sum_{a \in \mathcal{A}(c)} \sum_{b \in \mathcal{B}} \sum_{\ell_{1} \in \mathcal{L}(0,a,b)} \sum_{h \in \mathcal{H}(c,\ell)} x_{0}(c,a,\ell_{1},b,\ell,h)}\right].
\label{eqn:flowqueue0a}
\end{align}

\paragraph{\textbf{(S4) Initial ambulance supply at stations:}}
For each $a \in \mathcal{A}(c_{0})$, $b \in \mathcal{B}$,
\begin{equation}
\label{eqn:dispatchbases0}
\sum_{h \in \mathcal{H}(c_{0},\ell_{0})} x_{0}(c_{0},a,b,\ell_{0},h) \ \ \leq \ \ A_{0}(a,b).
\end{equation}

\paragraph{\textbf{(S5) Initial ambulance supply at locations:}}
For each $a \in \mathcal{A}(c_{0})$, $b \in \mathcal{B}$, $\ell_{1} \in \mathcal{L}(0,a,b)$,
\begin{equation}
\label{eqn:dispatchlocals0}
\sum_{h \in \mathcal{H}(c_{0},\ell_{0})} x_{0}(c_{0},a,\ell_{1},b,\ell_{0},h) \ \ \leq \ \ A_{0}(a, \ell_{1}, b).
\end{equation}

\paragraph{\textbf{(S6) Initial dispatch bound:}}
\begin{equation}
\displaystyle{\sum_{a \in \mathcal{A}(c_{0})} \sum_{b \in \mathcal{B}} \sum_{h \in \mathcal{H}(c_{0},\ell_{0})} x_{0}(c_{0},a,b,\ell_{0},h)} + \displaystyle{\sum_{a \in \mathcal{A}(c_{0})} \sum_{b \in \mathcal{B}} \sum_{\ell_{1} \in \mathcal{L}(0,a,b)} \sum_{h \in \mathcal{H}(c_{0},\ell_{0})} x_{0}(c_{0},a,\ell_{1},b,\ell_{0},h)} \ \ \le \ \ C_{0}(c_{0},\ell_{0})
\label{eqn:dispatchbound0a}
\end{equation}

\noindent
For $t = 1,\ldots,T$, the following constraints apply:

\paragraph{\textbf{(At1) Flow balance equations at the stations:}}
For each $t = 1,\ldots,T$, $a \in \mathcal{A}$, $b \in \mathcal{B}$,
\begin{align}
A_{t+1}(a,b) \ \ = \ \ & A_{t}(a,b) + \displaystyle{\sum_{\{h \in \mathcal{H} \; : \; L(t,a,h,b) = b\}} y_{t}(a,h,b)} \nonumber \\
& {} + \displaystyle{\sum_{\{\ell_{1} \in \mathcal{L}(t,a,b) \; : \; L(t,a,\ell_{1},b) = b\}} \left[A_{t}(a,\ell_{1},b) - \sum_{c \in \mathcal{C}(a)} \sum_{\ell \in \mathcal{L}} \sum_{h \in H(c,\ell)} x_{t}(c,a,\ell_{1},b,\ell,h)\right]} \nonumber \\
& {} - \displaystyle{\sum_{c \in \mathcal{C}(a)} \sum_{\ell \in \mathcal{L}} \sum_{h \in \mathcal{H}(c,\ell)} x_{t}(c,a,b,\ell,h)}.
\label{eqn:flowbaset}
\end{align}

\paragraph{\textbf{(At2) Flow balance equations at the hospitals:}}
For each $t = 1,\ldots,T$, $a \in \mathcal{A}$, $h \in \mathcal{H}$,
\begin{align}
& \displaystyle{\sum_{c \in \mathcal{C}(a)} \sum_{\ell \in \mathcal{L}} \sum_{h' \in \mathcal{H}(c,\ell)} x_{t}(c,a,h,\ell,h') + \sum_{b \in \mathcal{B}} y_{t}(a,h,b)} \nonumber \\
= \ \ & \displaystyle{\sum_{c \in \mathcal{C}(a)} \sum_{\ell_{1} \in \mathcal{L}} \sum_{\{\ell \in \mathcal{L} \; : \; h \in \mathcal{H}(c,\ell), \; \tau(0,c,a,\ell_{1},\ell,h) = t\}} A_{0}(c,a,\ell_{1},\ell,h)} \nonumber \\
& {} + \displaystyle{\sum_{c \in \mathcal{C}(a)} \sum_{\{\ell_{1} \in \mathcal{L} \; : \; \tau_{0}(c,a,\ell_{1},h) = t\}} A_{0}(c,a,\ell_{1},h)} \nonumber \\
& {} + \mathds{1}_{\{c_{0} \in \mathcal{C}(a), \; h \in \mathcal{H}(c_{0},\ell_{0})\}} \displaystyle{\sum_{b \in \mathcal{B}} \mathds{1}_{\{\tau(0,c_{0},a,b,\ell_{0},h) = t\}} x_{0}(c_{0},a,b,\ell_{0},h)} \nonumber \\
& {} + \displaystyle{\sum_{c \in \mathcal{C}(a)} \sum_{\{\ell \in \mathcal{L} \; : \; h \in \mathcal{H}(c,\ell)\}} \sum_{b \in \mathcal{B}} \sum_{\{t' \in \{1,\ldots,t-1\} \; : \; t' + \tau(t',c,a,b,\ell,h) = t\}} x_{t'}(c,a,b,\ell,h)} \nonumber \\
& {} + \displaystyle{\sum_{c \in \mathcal{C}(a)} \sum_{\{\ell \in \mathcal{L} \; : \; h \in \mathcal{H}(c,\ell)\}} \sum_{h' \in \mathcal{H}} \sum_{\{t' \in \{1,\ldots,t-1\} \; : \; t' + \tau(t',c,a,h',\ell,h) = t\}} x_{t'}(c,a,h',\ell,h)} \nonumber \\
& {} + \mathds{1}_{\{c_{0} \in \mathcal{C}(a), \; h \in \mathcal{H}(c_{0},\ell_{0})\}} \displaystyle{\sum_{b \in \mathcal{B}} \sum_{\ell_{1} \in \mathcal{L}(0,a,b)} \mathds{1}_{\{\tau(0,c_{0},a,\ell_{1},\ell_{0},h) = t\}} x_{0}(c_{0},a,\ell_{1},b,\ell_{0},h)} \nonumber \\
& {} + \displaystyle{\sum_{c \in \mathcal{C}(a)} \sum_{\{\ell \in \mathcal{L} \; : \; h \in \mathcal{H}(c,\ell)\}} \sum_{b \in \mathcal{B}} \sum_{\ell_{1} \in \mathcal{L}} \sum_{\{t' \in \{1,\ldots,t-1\} \; : \; \ell_{1} \in \mathcal{L}(t',a,b), \; t' + \tau(t',c,a,\ell_{1},\ell,h) = t\}} x_{t'}(c,a,\ell_{1},b,\ell,h)}.
\label{eqn:flowhosptcall0}
\end{align}

\paragraph{\textbf{(At3) Flow balance equations at the locations between hospitals and stations:}}
For each $t = 1,\ldots,T$, $a \in \mathcal{A}$, $b \in \mathcal{B}$, $\ell_{1} \in \mathcal{L}(t,a,b)$,
\begin{align}
A_{t+1}(a,\ell_{1},b) \ \ = \ \ & \displaystyle{\sum_{\{h \in \mathcal{H} \; : \; L(t,a,h,b) = \ell_{1}\}} y_{t}(a,h,b)} \nonumber \\
& {} + \displaystyle{\sum_{\{\ell'_{1} \in \mathcal{L}(t,a,b) \; : \; L(t,a,\ell'_{1},b) = \ell_{1}\}} \left[A_{t}(a,\ell'_{1},b) - \sum_{c \in \mathcal{C}(a)} \sum_{\ell \in \mathcal{L}} \sum_{h \in H(c,\ell)} x_{t}(c,a,\ell'_{1},b,\ell,h)\right]}.
\label{eqn:flowlocationt}
\end{align}

\paragraph{\textbf{(At4) Flow balance equations for the queues:}}
For each $t = 1,\ldots,T$, $c \in \mathcal{C}$, $\ell \in \mathcal{L}$,
\begin{align}
C_{t+1}(c,\ell) \ \ = \ \ & C_{t}(c,\ell) + \lambda(t,c,\ell) - \displaystyle{\sum_{a \in \mathcal{A}(c)} \sum_{b \in \mathcal{B}} \sum_{h \in \mathcal{H}(c,\ell)} x_{t}(c,a,b,\ell,h)} \nonumber \\
& {} - \displaystyle{\sum_{a \in \mathcal{A}(c)} \sum_{h' \in \mathcal{H}} \sum_{h \in \mathcal{H}(c,\ell)} x_{t}(c,a,h',\ell,h)} - \displaystyle{\sum_{a \in \mathcal{A}(c)} \sum_{b \in \mathcal{B}} \sum_{\ell_{1} \in \mathcal{L}(t,a,b)} \sum_{h \in H(c,\ell)} x_{t}(c,a,\ell_{1},b,\ell,h)}.
\label{eqn:flowqueuet}
\end{align}

\paragraph{\textbf{(At5) Station capacities:}}
For each $t = 1,\ldots,T+1$, $b \in \mathcal{B}$,
\begin{equation}
\label{eqn:capacity}
\sum_{a \in \mathcal{A}} A_{t}(a,b) \ \ \leq \ \ A_{\max}(b),
\end{equation}
where ${A}_{\max}(b)$ denotes the maximum number of ambulances that can park at station~$b$.

\paragraph{\textbf{(At6) Ambulance supply at stations:}}
For each $t = 1,\ldots,T$, $a \in \mathcal{A}$, $b \in \mathcal{B}$,
\begin{equation}
\label{eqn:dispatchbases}
\sum_{c \in \mathcal{C}(a)} \sum_{\ell \in \mathcal{L}} \sum_{h \in \mathcal{H}(c,\ell)} x_{t}(c,a,b,\ell,h) \ \ \leq \ \ A_{t}(a,b).
\end{equation}

\paragraph{\textbf{(At7) Ambulance supply at the locations between hospitals and stations:}}
For each $t = 1,\ldots,T$, $a \in \mathcal{A}$, $b \in \mathcal{B}$, $\ell_{1} \in \mathcal{L}(t,a,b)$,
\begin{equation}
\label{eqn:dispatchlocals}
\sum_{c \in \mathcal{C}(a)} \sum_{\ell \in \mathcal{L}} \sum_{h \in \mathcal{H}(c,\ell)} x_{t}(c,a,\ell_{1},b,\ell,h) \ \ \leq \ \ A_{t}(a, \ell_{1}, b).
\end{equation}

\paragraph{\textbf{(At8) Dispatch bounds:}}
For each $t = 1,\ldots,T$, $c \in \mathcal{C}$, $\ell \in \mathcal{L}$,
\begin{align}
& \displaystyle{\sum_{a \in \mathcal{A}(c)} \sum_{b \in \mathcal{B}} \sum_{h \in \mathcal{H}(c,\ell)} x_{t}(c,a,b,\ell,h)} + \displaystyle{\sum_{a \in \mathcal{A}(c)} \sum_{h' \in \mathcal{H}} \sum_{h \in \mathcal{H}(c,\ell)} x_{t}(c,a,h',\ell,h)} \nonumber \\
& {} + \displaystyle{\sum_{a \in \mathcal{A}(c)} \sum_{b \in \mathcal{B}} \sum_{\ell_{1} \in \mathcal{L}(t,a,b)} \sum_{h \in H(c,\ell)} x_{t}(c,a,\ell_{1},b,\ell,h)}
\ \ \le \ \ C_{t}(c,\ell) + \lambda(t,c,\ell)
\label{eqn:dispatchboundt}
\end{align}

\subsubsection{Constraints for the ambulance reassignment problem}

The following five sets of constraints apply to the first-stage variables for the ambulance reassignment problem only:

\paragraph{\textbf{(R1) Flow balance equations at the stations:}}
For each $a \in \mathcal{A}$, $b \in \mathcal{B}$,
\begin{equation}
A_{1}(a,b) \ \ = \ \ A_{0}(a,b) + \mathds{1}_{\{a = a_{0}, \, L(0,a_{0},h_{0},b) = b\}} y_{0}(a_{0},h_{0},b) + \displaystyle{\sum_{\{\ell_{1} \in \mathcal{L}(0,a,b) \; : \; L(0,a,\ell_{1},b) = b\}}} A_{0}(a,\ell_{1},b).
\label{eqn:flowbase0b}
\end{equation}

\paragraph{\textbf{(R2) Flow balance equation at hospital~$h_{0}$:}}
\begin{equation}
\displaystyle{\sum_{c \in \mathcal{C}(a_{0})} \sum_{\ell \in \mathcal{L}} \sum_{h \in \mathcal{H}(c,\ell)} x_{0}(c,a_{0},h_{0},\ell,h) + \sum_{b \in \mathcal{B}} y_{0}(a_{0},h_{0},b)} \ \ = \ \ 1.
\label{eqn:flowhosp0b}
\end{equation}

\paragraph{\textbf{(R3) Flow balance equations at the locations between hospitals and stations:}}
For each $a \in \mathcal{A}$, $b \in \mathcal{B}$, $\ell_{1} \in \mathcal{L}(0,a,b)$,
\begin{equation}
A_{1}(a,\ell_{1},b) \ \ = \ \ \displaystyle{\sum_{\{\ell'_{1} \in \mathcal{L}(0,a,b) \; : \; L(0,a,\ell'_{1},b) = \ell_{1}\}} A_{0}(a,\ell'_{1},b)} + \mathds{1}_{\{a = a_{0}, \, L(0,a_{0},h_{0},b) = \ell_{1}\}} y_{0}(a_{0},h_{0},b).
\label{eqn:flowlocation0b}
\end{equation}

\paragraph{\textbf{(R4) Flow balance equations for the queues:}}
For each $c \in \mathcal{C}$, $\ell \in \mathcal{L}$,
\begin{equation}
C_{1}(c,\ell) \ \ = \ \ C_{0}(c,\ell) - \mathds{1}_{\{c \in \mathcal{C}(a_{0})\}} \displaystyle{\sum_{h \in \mathcal{H}(c,\ell)} x_{0}(c,a_{0},h_{0},\ell,h)}.
\label{eqn:flowqueue0b}
\end{equation}

\paragraph{\textbf{(R5) Initial dispatch bounds:}}
For each $c \in \mathcal{C}(a_{0})$, $\ell \in \mathcal{L}$,
\begin{equation}
\displaystyle{\sum_{h \in \mathcal{H}(c,\ell)} x_{0}(c,a_{0},h_{0},\ell,h)} \ \ \le \ \ C_{0}(c,\ell).
\label{eqn:dispatchbound0b}
\end{equation}

\noindent
For $t=1,\ldots,T$, constraints (At1), (At3), (At4), (At5), (At6), (At7), and (At8) apply, and flow balance constraints (At2) at the hospitals are replaced by constraints (At9).

\paragraph{\textbf{(At9) Flow balance equations at the hospitals:}}
For each $t = 1,\ldots,T$, $a \in \mathcal{A}$, $h \in \mathcal{H}$,
\begin{align}
& \displaystyle{\sum_{c \in \mathcal{C}(a)} \sum_{\ell \in \mathcal{L}} \sum_{h' \in \mathcal{H}(c,\ell)} x_{t}(c,a,h,\ell,h')} + \displaystyle{\sum_{b \in \mathcal{B}} y_{t}(a,h,b)} \nonumber \\
= \ \ & \displaystyle{\sum_{c \in \mathcal{C}(a)} \sum_{\ell_{1} \in \mathcal{L}} \sum_{\{\ell \in \mathcal{L} \; : \; h \in \mathcal{H}(c,\ell), \; \tau(0,c,a,\ell_{1},\ell,h) = t\}} A_{0}(c,a,\ell_{1},\ell,h)} \nonumber \\
& {} + \displaystyle{\sum_{c \in \mathcal{C}(a)} \sum_{\{\ell_{1} \in \mathcal{L} \; : \; \tau_{0}(c,a,\ell_{1},h) = t\}} A_{0}(c,a,\ell_{1},h)} \nonumber \\
& {} + \displaystyle{\sum_{c \in \mathcal{C}(a)} \sum_{\{\ell \in \mathcal{L} \; : \; h \in \mathcal{H}(c,\ell)\}} \sum_{b \in \mathcal{B}} \sum_{\{t' \in \{1,\ldots,t-1\} \; : \; t' + \tau(t',c,a,b,\ell,h) = t\}} x_{t'}(c,a,b,\ell,h)} \nonumber \\
& {} + \mathds{1}_{\{a = a_{0}\}} \displaystyle{\sum_{c \in \mathcal{C}(a)} \sum_{\{\ell \in \mathcal{L} \; : \; h \in \mathcal{H}(c,\ell)\}} \mathds{1}_{\{\tau(0,c,a_{0},h_{0},\ell,h) = t\}} x_{0}(c,a_{0},h_{0},\ell,h)} \nonumber \\
& {} + \displaystyle{\sum_{c \in \mathcal{C}(a)} \sum_{\{\ell \in \mathcal{L} \; : \; h \in \mathcal{H}(c,\ell)\}} \sum_{h' \in \mathcal{H}} \sum_{\{t' \in \{1,\ldots,t-1\} \; : \; t' + \tau(t',c,a,h',\ell,h) = t\}} x_{t'}(c,a,h',\ell,h)} \nonumber \\
& {} + \displaystyle{\sum_{c \in \mathcal{C}(a)} \sum_{\{\ell \in \mathcal{L} \; : \; h \in \mathcal{H}(c,\ell)\}} \sum_{b \in \mathcal{B}} \sum_{\ell_{1} \in \mathcal{L}} \sum_{\{t' \in \{1,\ldots,t-1\} \; : \; \ell_{1} \in \mathcal{L}(t',a,b), \; t' + \tau(t',c,a,\ell_{1},\ell,h) = t\}} x_{t'}(c,a,\ell_{1},b,\ell,h)}.
\label{eqn:flowhosptamb0}
\end{align}

\subsubsection{Objective function for the ambulance selection problem}

Let $f_{t}(c,a,b,\ell,h)$ denote the cost per emergency if ambulance~$a$ is dispatched at time~$t$ from station~$b$ to serve an emergency of type~$c$ at location~$\ell$, and transport patient(s) to hospital~$h$; this includes a penalty for the waiting time and other costs.
Similarly, let $f_{t}(c,a,h',\ell,h)$ denote the cost per emergency if ambulance~$a$ is dispatched at time~$t$ from hospital~$h'$ to serve an emergency of type~$c$ at location~$\ell$, and transport patient(s) to hospital~$h$;
let $f_{t}(c,a,\ell_{1},b,\ell,h)$ denote the cost per emergency if ambulance~$a$ at location~$\ell_{1}$ at time~$t$ traveling toward station~$b$ is dispatched to serve an emergency of type~$c$ that arrived at location~$\ell$ at time~$t$, and transport patient(s) to hospital~$h$;
let $f_{t}(a,h,b)$ denote the cost if ambulance~$a$ is dispatched at time~$t$ from hospital~$h$ to station~$b$;
let $g_{t}(c,\ell)$ denote the penalty per emergency of type~$c$ waiting in queue at location~$\ell$ at the beginning of time~$t$;
let $g_{t}(a,b)$ denote the cost if ambulance~$a$ is at station~$b$ at the beginning of time~$t$;
and let $g_{t}(a,\ell_{1},b)$ denote the cost if ambulance~$a$ moves from location~$\ell_{1}$ to location~$L(t,a,\ell_{1},b)$ during time~$t$.

For the ambulance selection problem, the objective is to minimize
\begin{align}
& \displaystyle{\sum_{a \in \mathcal{A}(c_{0})} \sum_{b \in \mathcal{B}} \sum_{h \in \mathcal{H}(c_{0},\ell_{0})} \left[f_{0}(c_{0},a,b,\ell_{0},h) x_{0}(c_{0},a,b,\ell_{0},h) + \sum_{\ell_{1} \in \mathcal{L}(0,a,b)} f_{0}(c_{0},a,\ell_{1},b,\ell_{0},h) x_{0}(c_{0},a,\ell_{1},b,\ell_{0},h)\right]} \nonumber \\
& {} + \displaystyle{\sum_{t=1}^{T} \sum_{c \in \mathcal{C}} \sum_{a \in \mathcal{A}(c)} \sum_{b \in \mathcal{B}} \sum_{\ell \in \mathcal{L}} \sum_{h \in \mathcal{H}(c,\ell)} \left[f_{t}(c,a,b,\ell,h) x_{t}(c,a,b,\ell,h) + \sum_{\ell_{1} \in \mathcal{L}(t,a,b)} f_{t}(c,a,\ell_{1},b,\ell,h) x_{t}(c,a,\ell_{1},b,\ell,h)\right]} \nonumber \\
& {} + \displaystyle{\sum_{t=1}^{T} \sum_{c \in \mathcal{C}} \sum_{a \in \mathcal{A}(c)} \sum_{\ell \in \mathcal{L}} \sum_{h \in \mathcal{H}(c,\ell)} \sum_{h' \in \mathcal{H}} f_{t}(c,a,h',\ell,h) x_{t}(c,a,h',\ell,h) + \sum_{t=1}^{T} \sum_{a \in \mathcal{A}} \sum_{b \in \mathcal{B}} \sum_{h \in \mathcal{H}} f_{t}(a,h,b) y_{t}(a,h,b)} \nonumber \\
& {} + \displaystyle{\sum_{t=0}^{T} \left[\sum_{c \in \mathcal{C}} \sum_{\ell \in \mathcal{L}} g_{t+1}(c,\ell) C_{t+1}(c,\ell) + \sum_{a \in \mathcal{A}} \sum_{b \in \mathcal{B}} \left\{g_{t+1}(a,b) A_{t+1}(a,b) + \sum_{\ell_{1} \in \mathcal{L}(t+1,a,b)} g_{t+1}(a,\ell_{1},b) A_{t+1}(a,\ell_{1},b)\right\}\right]}.
\label{eqn:objectivecall0}
\end{align}

\subsubsection{Objective function for the ambulance reassignment problem}

Using the notation of the previous sections, for the ambulance reassignment problem, the objective is to minimize
\begin{align}
& \displaystyle{\sum_{c \in \mathcal{C}(a_{0})} \sum_{\ell \in \mathcal{L}} \sum_{h \in \mathcal{H}\mathcal{H}(c,\ell)} f_{0}(c,a_{0},h_{0},\ell,h) x_{0}(c,a_{0},h_{0},\ell,h) + \sum_{b \in \mathcal{B}} f_{0}(a_{0},h_{0},b) y_{0}(a_{0},h_{0},b)} \nonumber \\
& \displaystyle{\sum_{t=1}^{T} \sum_{c \in \mathcal{C}} \sum_{a \in \mathcal{A}(c)} \sum_{b \in \mathcal{B}} \sum_{\ell \in \mathcal{L}} \sum_{h \in \mathcal{H}(c,\ell)} \left[f_{t}(c,a,b,\ell,h) x_{t}(c,a,b,\ell,h) + \sum_{\ell_{1} \in \mathcal{L}(t,a,b)} f_{t}(c,a,\ell_{1},b,\ell,h) x_{t}(c,a,\ell_{1},b,\ell,h)\right]} \nonumber \\
& {} + \displaystyle{\sum_{t=1}^{T} \left[\sum_{c\in \mathcal{C}} \sum_{a \in \mathcal{A}(c)} \sum_{h' \in \mathcal{H}} \sum_{\ell \in \mathcal{L}} \sum_{h \in \mathcal{H}(c,\ell)} f_{t}(c,a,h',\ell,h) x_{t}(c,a,h',\ell,h) + \sum_{a \in \mathcal{A}} \sum_{h \in \mathcal{H}} \sum_{b \in \mathcal{B}} f_{t}(a,h,b) y_{t}(a,h,b)\right]} \nonumber \\
& {} + \displaystyle{\sum_{t=0}^{T} \left[\sum_{c \in \mathcal{C}} \sum_{\ell \in \mathcal{L}} g_{t+1}(c,\ell) C_{t+1}(c,\ell) + \sum_{a \in \mathcal{A}} \sum_{b \in \mathcal{B}} \left\{g_{t+1}(a,b) A_{t+1}(a,b) + \sum_{\ell_{1} \in \mathcal{L}(t+1,a,b)} g_{t+1}(a,\ell_{1},b) A_{t+1}(a,\ell_{1},b)\right\}\right]}.
\label{eqn:objectiveamb0}
\end{align}

\subsection{Deterministic Itinerary-based Model}
\label{sec:detitinerarymodel}

\subsubsection{Problem input parameters}
\label{sec:parameters}

At time~$0$ the current state of each ambulance $a \in \mathcal{A}$ is known.
Let $\mathcal{E}$ denote the set of emergencies from time~$0$ onward, including the emergencies waiting in queue at time~$0$, any call that has just arrived at time~$0$, and future emergencies that will arrive after time~$0$.
For each emergency $e \in \mathcal{E}$, let $t_{1}(e)$ denote the arrival time of emergency~$e$, let $c(e) \in \mathcal{C}$ denote the type of emergency~$e$, and let $\ell(e) \in \mathcal{L}$ denote the location of emergency~$e$.
Also, for each ambulance $a \in \mathcal{A}$, let $\mathcal{E}(a) \defi \{e \in \mathcal{E} \, : \, c(e) \in \mathcal{C}(a)\}$ denote the set of emergencies that can be handled by ambulance~$a$.

Next we describe an itinerary for an ambulance $a \in \mathcal{A}$.
Each itinerary starts at time~$0$.
If the ambulance was dispatched to an emergency before time~$0$ and it has not yet completed the task for the emergency, then it can be determined when and where the ambulance will become available for its next dispatch.
Otherwise, the ambulance is available to be dispatched at time~$0$.
Thereafter an itinerary consists of a sequence of dispatches to either emergencies $e \in \mathcal{E}(a)$ or to ambulance stations $b \in \mathcal{B}$.
Given a sequence of dispatches up to and including a particular emergency $e \in \mathcal{E}(a)$ in an itinerary for a particular ambulance $a \in \mathcal{A}$, it can be determined at what time~$t_{2}(e)$ and where ambulance~$a$ will become available after completing its task for emergency~$e$.
Then ambulance~$a$ must either be dispatched to an emergency that has already arrived and that can be handled by the ambulance, that is, an emergency in $\{e' \in \mathcal{E}(a) \, : \, t_{1}(e') < t_{2}(e)\}$, or the ambulance must be dispatched to an ambulance station $b \in \mathcal{B}$.
Note that although all emergencies $\mathcal{E}$ are included in the input of the deterministic formulation, an ambulance cannot be dispatched to an emergency before the emergency arrives.
The ambulance station to which an ambulance is dispatched can be chosen with knowledge of future calls, and thus the deterministic formulation still provides an optimistic bound of the objective value of the stochastic problem.
If an ambulance~$a$ is dispatched to an ambulance station~$b$ at time~$t_{2}(e)$, then thereafter the ambulance can be dispatched to any emergency $e' \in \mathcal{E}(a)$ that has not been served yet.
If $t_{1}(e')$ is less than the arrival time of ambulance~$a$ to station~$b$, then the ambulance can be dispatched to emergency~$e'$ while still on the way to~$b$, and the arrival time of ambulance~$a$ at emergency~$e'$ can be calculated accordingly.
If $t_{1}(e')$ is greater than the arrival time of ambulance~$a$ to station~$b$, then the arrival time of ambulance~$a$ at emergency~$e'$ is no earlier than $t_{1}(e')$ plus the travel time from~$b$ to $\ell(e')$.

For each ambulance $a \in \mathcal{A}$, let $\mathcal{I}(a)$ denote the set of such itineraries for ambulance~$a$ over the time horizon $[0,T]$, and let $\mathcal{I} \defi \cup_{a \in \mathcal{A}} \mathcal{I}(a)$ denote the set of feasible itineraries.
An ambulance can serve approximately one emergency per hour, and thus if the time horizon is a few hours long, the number $\left|\mathcal{I}(a)\right|$ of itineraries for an ambulance is not very large.
For any ambulance $a \in \mathcal{A}$ and emergency $e \in \mathcal{E}(a)$, let $\mathcal{I}_{e}(a) \subset \mathcal{I}(a)$ denote the set of all itineraries of ambulance~$a$ that serve emergency~$e$.
For each ambulance $a \in \mathcal{A}$, emergency $e \in \mathcal{E}(a)$, and itinerary $i \in \mathcal{I}_{e}(a)$, the response time for emergency~$e$ can be determined, and the quality-of-care obtained by matching ambulance~$a$ with emergency~$e$ can be determined.
That enables determination of a cost $f_{i}$ for each itinerary $i \in \mathcal{I}$.
In addition, for each emergency $e \in \mathcal{E}$, there is a penalty $g_{e}$ if emergency~$e$ is not served by any ambulance during the time horizon.

\subsubsection{Itinerary-based optimization problem}

For each itinerary $i \in \mathcal{I}$, let decision variable $z_{i} = 1$ if itinerary $i$ is chosen, and $z_{i} = 0$ otherwise.
Then the itinerary-based optimization problem is as follows:
\begin{align}
\label{eqn:itinerary opt problem objective}
\min \ \ & \sum_{i \in \mathcal{I}} f_{i} z_{i} + \sum_{e \in \mathcal{E}} g_{e} \left(1 - \sum_{a \in \mathcal{A}(c(e))} \sum_{i \in \mathcal{I}_{e}(a)} z_{i}\right) \\
\label{eqn:itinerary opt problem each ambulance itinerary}
\mbox{s.t.} \ \ & \sum_{i \in \mathcal{I}(a)} z_{i} \ \ = \ \ 1 \qquad \forall \ a \in \mathcal{A} \\
& \sum_{a \in \mathcal{A}(c(e))} \sum_{i \in \mathcal{I}_{e}(a)} z_{i} \ \ \le \ \ 1 \qquad \forall \ e \in \mathcal{E}
\label{eqn:itinerary opt problem each emergency one visit}
\end{align}

\ignore{
Let $x^*$ denote an optimal solution of problem~\eqref{eqn:itinerary opt problem}.
If the problem on-hand is an ambulance selection problem, then let $e_{0} \in \mathcal{E}$ denote the emergency that has just arrived at time~$0$.
If $x^*_{i} = 1$ for some itinerary $i \in \mathcal{I}_{e_{0}}(a)$, and itinerary~$i$ dispatches ambulance~$a$ to emergency~$e_{0}$ at time~$0$, then the corresponding policy implements such a dispatch decision at time~$0$.
Otherwise the policy places emergency~$e_{0}$ in queue, and the ambulance to be dispatched to the emergency will be determined by solving an ambulance reassignment problem later.
If the problem on-hand is an ambulance reassignment problem, then let $a_{0} \in \mathcal{A}$ denote the ambulance that has just become available at time~$0$.
Consider the unique itinerary $i \in \mathcal{I}(a_{0})$ such that $x^*_{i} = 1$.
Then itinerary~$i$ must start with a dispatch of ambulance~$a_{0}$ at time~$0$, either to an emergency that is in queue at time~$0$, or to an ambulance station.
The policy implements that dispatch immediately.
}

\subsection{Stochastic Model}
\label{sec:smodel}

In this section we describe a stochastic model of ambulance dispatch operations.
The stochastic model is used to generate scenarios for the two-stage optimization problem, and is used in a simulation to test the performance of different ambulance dispatch policies.

The same space discretization is used for the optimization problem and the simulation, but the time discretization applies to the arc-based optimization problem only.
The random variables are the sequences of emergency calls of each type at each location, the ambulance travel times, and the emergency service times.
Emergency calls arrive in continuous time according to an exogenous stochastic process, that is, we assume that the emergency arrivals do not depend on the dispatch decisions.
We assume that the stochastic process has the property that with probability~$1$, at most one emergency call arrives at a continuous time point.
Let $\omega(\tau)$ denote a random sequence of emergency calls, travel time random variables, and service time random variables, during time interval $(\tau,T')$, where $T' \gg T$ denotes a specified simulation time horizon, and let $\xi(\tau)$ denote a random sequence of emergency calls, travel time random variables, and service time random variables, during time interval $(\tau,\tau+T)$.
For any time~$\tau$, let $s(\tau)$ denote the state of the process at time~$\tau$.
That is, $s(\tau)$ contains information about the location and current assignment of each ambulance, and the calls in queue at time~$\tau$, including the data of a call, if any, that has just arrived.
Let $\pi$ denote any deterministic policy that takes the state $s(\tau)$ at any time~$\tau$ as input, and specifies either an ambulance selection decision or an ambulance reassignment decision $\pi(s(\tau))$, as appropriate.
The simulation generates a sample path $\omega(0)$ according to the specified stochastic process, and keeps track of the state and performance metrics of the process.
When an ambulance selection decision or an ambulance reassignment decision has to be made at time~$\tau$, the simulation provides $s(\tau)$ as input to $\pi$, receives the decision $\pi(s(\tau))$ as output from $\pi$, and updates the state of the process accordingly.

Next we describe how the deterministic optimization problems specified in the previous sections are used to compute a policy $\pi$.
First, suppose that $\pi$ is called at time~$\tau$ with a request for an ambulance selection decision, and with state $s(\tau)$ provided as input.
Then time~$\tau$ in the simulation corresponds to time $t=0$ in the current optimization problem.
Also, $N$ independent and identically distributed second-stage scenarios, denoted by $\xi_{1}(\tau),\dots,\xi_{N}(\tau)$, is generated according to the specified stochastic process.
Note that scenarios $\xi_{1}(\tau),\dots,\xi_{N}(\tau)$ are independent of the sample path $\omega(\tau)$ used to evaluate the policies.
Also note that all the input of the optimization problem either does not change, or can be derived from the state $s(\tau)$ provided as input by the simulation and the generated scenarios $\xi_{1}(\tau),\dots,\xi_{N}(\tau)$.
For example, the set $\mathcal{C}$ of emergency types does not change, the type $c_{0}$ of the call at $t=0$ is part of the state $s(\tau)$, and the sequence of emergency calls of each type~$c$ at each location~$\ell$ during time interval $(\tau,\tau+T)$ under scenario~$n$ can be derived from $\xi_{n}(\tau)$.

Let $x_{0} \defi \big(x_{0}(c_{0},a,b,\ell_{0},h), x_{0}(c_{0},a,\ell_{1},b,\ell_{0},h), a \in \mathcal{A}(c_{0}), b \in \mathcal{B}, \ell_{1} \in \mathcal{L}(0,a,b), h \in \mathcal{H}(c_{0},\ell_{0})\big)$ denote the first-stage decision variables of the ambulance selection problem.
Then, let
\begin{align*}
F(x_{0}) \ \ \defi \ \ & \sum_{a \in \mathcal{A}(c_{0})} \sum_{h \in \mathcal{H}(c_{0},\ell_{0})} \sum_{b \in \mathcal{B}} \bigg[f_{0}(c_{0},a,b,\ell_{0},h) x_{0}(c_{0},a,b,\ell_{0},h) \\
& \hspace{35mm} {} + \sum_{\ell_{1} \in \mathcal{L}(0,a,b)} f_{0}(c_{0},a,\ell_{1},b,\ell_{0},h) x_{0}(c_{0},a,\ell_{1},b,\ell_{0},h)\bigg]
\end{align*}
denote the first-stage part of the ambulance selection problem objective.

For the arc-based model, and for each scenario~$n$, let $(x_{n},y_{n}) \defi \big(x_{n,t}(c,a,b,\ell,h), x_{n,t}(c,a,\ell_{1},b,\ell,h), x_{n,t}(c,a,h',\ell,h), y_{n,t}(a,h',b), C_{n,t}(c,\ell), A_{n,t}(a,b), A_{n,t}(a,\ell_{1},b), t = 1,\dots,T, c \in \mathcal{C}, a \in \mathcal{A}(c), b \in \mathcal{B}, \ell \in \mathcal{L}, \ell_{1} \in \mathcal{L}(0,a,b), h \in \mathcal{H}(c,\ell), h' \in \mathcal{H}\big)$ denote the second-stage decision variables for scenario~$n$.
Then, for each scenario~$n$, let $G(s(\tau),\xi_{n}(\tau),x_{n},y_{n})$ denote the second-stage part of the objective, that is the same for the ambulance selection problem and the ambulance reassignment problem.
For each first-stage decision $x_{0}$ and each scenario~$n$, the second-stage of the ambulance selection problem is
\begin{align}
Q(s(\tau),\xi_{n}(\tau),x_{0}) \ \ \defi \ \ \displaystyle{\min_{x_{n},y_{n}}} \ \ & G(s(\tau),\xi_{n}(\tau),x_{n},y_{n}) \nonumber \\
\mbox{s.t.} \ \ & (At1),(At2),(At3),(At4),(At5),(At6),(At7),(At8).
\label{eqn:arc second-stage select}
\end{align}
Policy~$\pi$ then chooses an ambulance selection decision
\begin{align}
\pi(s(\tau)) \ \ \in \ \ \argmin_{x_{0}} \ \ & F(x_{0}) + \frac{1}{N} \sum_{n=1}^{N} Q\big(s(\tau),\xi_{n}(\tau),x_{0}\big) \nonumber \\
\mbox{s.t.} \ \ & (S1),(S2),(S3),(S4),(S5),(S6).
\label{eqn:selection policy}
\end{align}
Similarly, suppose that $\pi$ is called at time~$\tau$ with a request for an ambulance reassignment decision, and with state $s(\tau)$ provided as input.
Let $\tilde{x}_{0} \defi \big(x_{0}(c,a_{0},h_{0},\ell,h), y_{t}(a_{0},h_{0},b), c \in \mathcal{C}, b \in \mathcal{B}, \ell \in \mathcal{L}, h \in \mathcal{H}(c,\ell)\big)$ denote the first-stage decision variables.
Then, let
\begin{align*}
\tilde{F}(\tilde{x}_{0}) \ \ \defi \ \ & \sum_{c\in \mathcal{C}} \sum_{\ell \in \mathcal{L}} \sum_{h \in \mathcal{H}(c,\ell)} f_{0}(c,a_{0},h_{0},\ell,h) x_{0}(c,a_{0},h_{0},\ell,h) + \sum_{b \in \mathcal{B}} f_{0}(a_{0},h_{0},b) y_{0}(a_{0},h_{0},b)
\end{align*}
denote the first-stage part of the ambulance reassignment problem objective.
For each first-stage decision $\tilde{x}_{0}$ and each scenario~$n$, the second-stage of the ambulance reassignment problem is
\begin{align}
\tilde{Q}(s(\tau),\xi_{n}(\tau),\tilde{x}_{0}) \ \ \defi \ \ \displaystyle{\min_{x_{n},y_{n}}} \ \ & G(s(\tau),\xi_{n}(\tau),x_{n},y_{n}) \nonumber \\
\mbox{s.t.} \ \ & (At1),(At3),(At4),(At5),(At6),(At7),(At8),(At9).
\label{eqn:arc second-stage reassign}
\end{align}
Policy~$\pi$ then chooses an ambulance reassignment decision
\begin{align}
\pi(s(\tau)) \ \ \in \ \ \argmin_{\tilde{x}_{0}} \ \ & \tilde{F}(\tilde{x}_{0}) + \frac{1}{N} \sum_{n=1}^{N} \tilde{Q}(s(\tau),\xi_{n}(\tau),\tilde{x}_{0}) \nonumber \\
\mbox{s.t.} \ \ & (R1),(R2),(R3),(R4),(R5).
\label{eqn:reassignment policy}
\end{align}

For the itinerary-based model, and for each scenario~$n$, let $z_{n} \defi \big(z_{n,i}, i \in \mathcal{I}\big)$ denote the second-stage decision variables for scenario~$n$.
Then, for each scenario~$n$, let $G(s(\tau),\xi_{n}(\tau),z_{n})$ denote the second-stage part of the objective, that is the same for the ambulance selection problem and the ambulance reassignment problem.
For each first-stage decision $x_{0}$ and each scenario~$n$, the second-stage of the ambulance selection problem is
\begin{align}
Q(s(\tau),\xi_{n}(\tau),x_{0}) \ \ \defi \ \ \displaystyle{\min_{z_{n}}} \ \ & G(s(\tau),\xi_{n}(\tau),z_{n}) \nonumber \\
\mbox{s.t.} \ \ & \eqref{eqn:itinerary opt problem each ambulance itinerary}, \eqref{eqn:itinerary opt problem each emergency one visit}.
\label{eqn:itinerary second-stage select}
\end{align}
Policy~$\pi$ then chooses an ambulance selection decision as given in~\eqref{eqn:selection policy}.

For each first-stage decision $\tilde{x}_{0}$ and each scenario~$n$, the second-stage of the ambulance reassignment problem is
\begin{align}
\tilde{Q}(s(\tau),\xi_{n}(\tau),\tilde{x}_{0}) \ \ \defi \ \ \displaystyle{\min_{z_{n}}} \ \ & G(s(\tau),\xi_{n}(\tau),z_{n}) \nonumber \\
\mbox{s.t.} \ \ & \eqref{eqn:itinerary opt problem each ambulance itinerary}, \eqref{eqn:itinerary opt problem each emergency one visit}.
\label{eqn:itinerary second-stage reassign}
\end{align}
Policy~$\pi$ then chooses an ambulance reassignment decision as given in~\eqref{eqn:reassignment policy}.

\section{Column Generation}
\label{sec:sol}

We propose a column generation algorithm to solve the optimization problems~\eqref{eqn:arc second-stage select} and~\eqref{eqn:arc second-stage reassign}.
Problems~\eqref{eqn:arc second-stage select} and~\eqref{eqn:arc second-stage reassign} are large linear programs that have to be solved fast in practice.
In addition, in both problems the number of decision variables is large relative to the number of constraints, so that in a basic feasible solution most decision variables have value zero.
This motivates the column generation algorithm proposed in this section to solve these problems.
We briefly describe the algorithm for problem~\eqref{eqn:arc second-stage select}.
A similar algorithm is used for problem~\eqref{eqn:arc second-stage reassign}.
Additional details are given in the Electronic Companion, Section~\ref{sec:column generation details}.

The column generation algorithm starts with an initial feasible solution for the continuous relaxation of problem~\eqref{eqn:arc second-stage select}.
Then an optimal primal-dual pair is computed for the restriction of the continuous relaxation of problem~\eqref{eqn:arc second-stage select} with only the decision variables that are nonzero in the initial feasible solution.
If the optimal dual solution of the restricted problem is feasible for the dual of the continuous relaxation of problem~\eqref{eqn:arc second-stage select}, then the primal-dual pair is optimal for the continuous relaxation of problem~\eqref{eqn:arc second-stage select}.
Otherwise, there exists a primal decision variable with negative reduced cost, and variables with negative reduced cost are added to the restricted problem to obtain the next restricted problem.
Every iteration adds a number of primal decision variables with negative reduced cost until the optimal dual solution of the restriction of the continuous relaxation of problem~\eqref{eqn:arc second-stage select} is feasible for the dual of the continuous relaxation of problem~\eqref{eqn:arc second-stage select}.
The Electronic Companion, Section~\ref{sec:column generation details}, gives an explanation how to find primal decision variables with negative reduced cost.

\ignore{
\subsection{Stochastic setting}

For the dispatch of ambulances
under uncertainty in the process
of calls, we have to solve in a rolling horizon approach sequences
of two-stage stochastic programs
as explained in Section \ref{sec:smodel}.
Since in our context there is
a small number of first stage
decisions, we solve each two-stage
problem enumerating all candidate
first stage decisions and computing
the second stage cost for every such candidate first stage decision.
The computation of the second
stage cost on a given scenario
is obtained solving a deterministic
problem (with calls given by the
calls of the corresponding scenario) which is either
a continuous relaxation of 
an ambulance selection problem
or a continuous relaxation of a
 ambulance reassignment problem.
We either solve these
problems directly or, if these
problems are of too large size and the assumptions (A1), (A2), (A3), (A4), (B1), (B2), (C1), (C2) of the previous section  apply,
using one of the column generation
Algorithms \ref{alg:columngeneration1},
\ref{alg:columngeneration2},
\ref{alg:columngeneration3},
\ref{alg:columngeneration4}, we have just described.
}

\section{Case Study}
\label{sec:num}

For the case study, an EMS system similar to Rio de Janeiro EMS (called the SAMU in what follows) was simulated.
All algorithms were implemented using C++14 and run on a computer with a Ryzen 5 2600 processor with 16GB of RAM memory, and a Ubuntu 22.04 OS.

\begin{figure}
\centering
\includegraphics[scale=0.6]{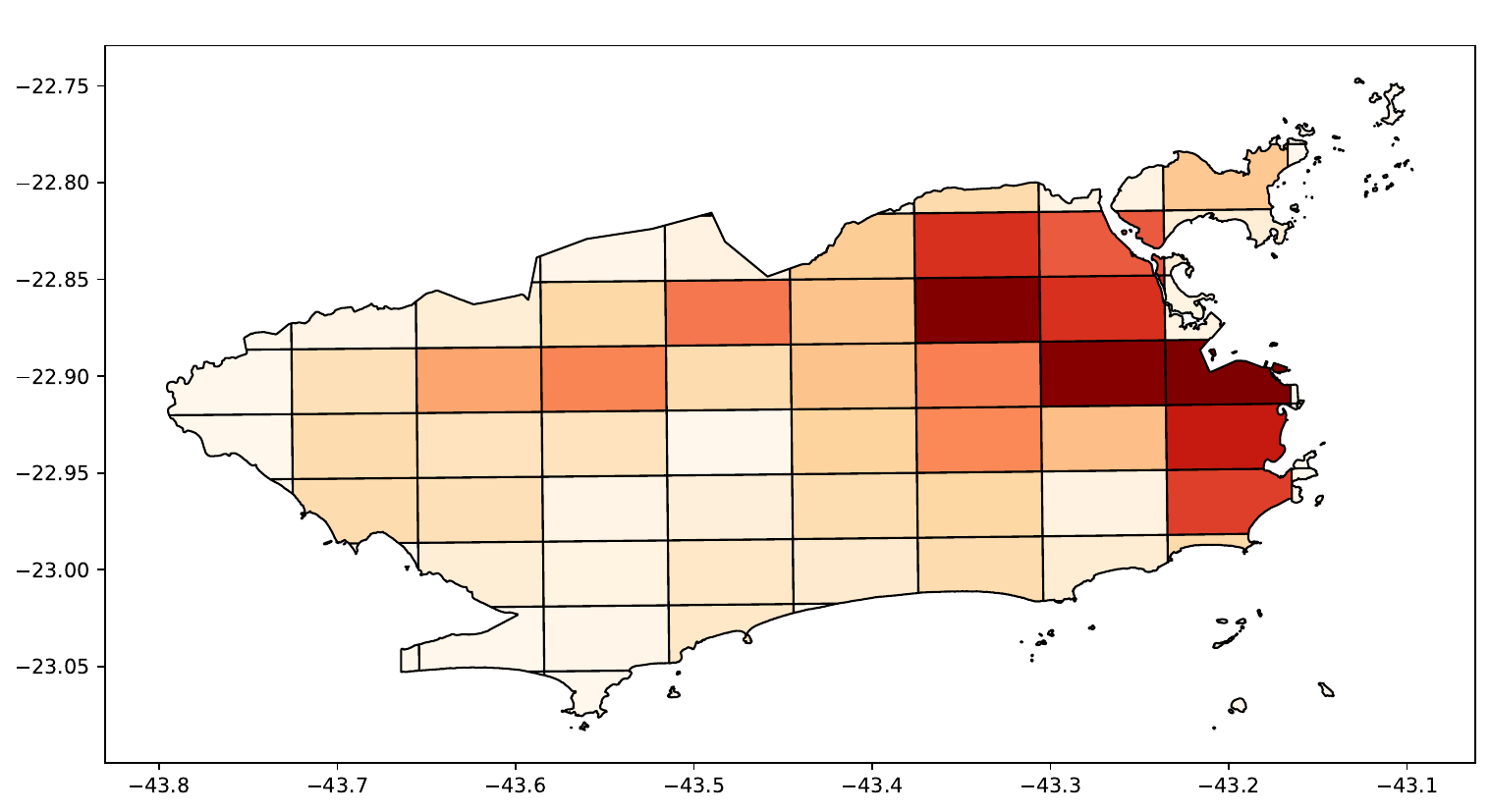}
\caption{Discretization of the service region of the Rio de Janeiro SAMU into $10 \times 10$ rectangles, and a heatmap of the mean intensity of emergency calls to this SAMU for the period January 2016--February 2018.\label{figuresheat}}
\end{figure}

Some features of the SAMU were simplified in the comparisons to enable use of the heuristics in the literature.
Ambulances/crews with two levels of capabilities (BLS and ALS) were simulated.
Emergencies were grouped into four types, consistent with the MPDS A/B/C/D classification:
\begin{itemize}
\item
Type~A: Low time urgency, no ambulance/crew capability preference;
\item
Type~B: High time urgency, no ambulance/crew capability preference;
\item
Type~C: Low time urgency, high ambulance/crew capability preferred;
\item
Type~D: High time urgency, high ambulance/crew capability preferred.
\end{itemize}
A rectangle that contains the service region of the SAMU was discretized into $10 \times 10$ rectangles.
Of these $100$~rectangles, $76$ have an intersection with the SAMU service region, and are shown in Figure~\ref{figuresheat}.
The SAMU has $10$~associated hospitals.
The number of ambulances and ambulance stations was varied from $10$ to $30$, with the number of ambulances always equal to the number of ambulance stations.
All patients were transported to the closest hospital.
For the arc-based model, time was discretized into $30$~minute intervals.

Emergency call data for the SAMU for the period January 2016--February 2018 were used to calibrate the call arrival rates $\lambda(t,c,\ell)$ for each $30$~minute period~$t$ of each day of the week, for each emergency type~$c$, and for each location~$\ell$, of a Poisson process, via maximizing a regularized likelihood function (see \cite{laspatedpaper,laspatedmanual,websiteambrouting24, ourheuristic241}).
Figure~\ref{figuresheat} also shows a heatmap of the total estimated rates of the Poisson process (summed over the emergency types and time periods) for each location~$\ell$.
For the numerical tests, we considered a time interval with high call arrival rates (Fridays between 18h00 and 20h00), and we simulated the call arrivals and dispatch decisions during this time interval.

Random travel times and service times were simulated as follows.
For each origin-destination pair $(a,b)$, the free-flow travel time $t_{f}(a,b)$ from~$a$ to~$b$, to travel the distance from~$a$ to~$b$ at a free-flow speed of 60km/h, was specified as input.
Given a uniform random number $u \in [0,1)$, the random travel time from~$a$ to~$b$ was given by $t(a,b,u) = t_{f}(a,b) - t_{f}(a,b) \ln(1-u)$, that is, the random travel time was equal to the free-flow travel time plus an exponentially distributed delay with mean equal to the free-flow travel time.
Let $c = 0$ for type~A emergencies, $c = 1$ for type~B emergencies, $c = 2$ for type~C emergencies, and $c = 3$ for type~D emergencies.
Then, given a uniform random number $u \in [0,1)$, the random on-site service time in seconds was given by $t_{\mathrm{scene}}(c,u) = 600 - 150 c \ln(1-u)$, that is, the random on-site service time was equal to $600$ seconds plus an exponentially distributed delay with mean equal to $150 c$ seconds.
Similarly, the random service time at a hospital was given by $t_{\mathrm{hosp}}(c,u) = 600 + \min\{600, - 150 c \ln(1-u)\}$.

\subsection{Comparison of Proposed Policies with Heuristics from the Literature}
\label{sec:policy comparison}

The performance of the optimization-based policies for ambulance dispatch decisions described in Section~\ref{sec:optmodel} was compared with the performance of policies proposed in the literature using simulation.
For each combination of number of ambulances and ambulance stations, and policy, we simulated $4$~replications of the time period Fridays between 18h00 and 20h00.
For the arc-based policy and the itinerary-based policy, we used $15$~second-stage scenarios (independent from the simulated scenario) with a time horizon length of $T = 2$~hours for each second-stage scenario.

The numerical results are summarized in Tables~\ref{tbl:enum_vs_all_random_nro}, \ref{tbl:enum_vs_all_actual_random_nro}, \ref{tbl:arc_vs_all_random_nro}, and~\ref{tbl:arc_vs_all_actual_random_nro}.
The tables have the same format: the first column gives the number of ambulances and ambulance stations, and the second column gives the average objective value (specified below) for the reference policy (the itinerary-based policy in Tables~\ref{tbl:enum_vs_all_random_nro} and~\ref{tbl:enum_vs_all_actual_random_nro}, and the arc-based policy in Tables~\ref{tbl:arc_vs_all_random_nro} and~\ref{tbl:arc_vs_all_actual_random_nro}).
The remaining columns give the results for other policies: ``Arc'' refers to the arc-based policy, ``Itinerary'' refers to the itinerary-based policy, ``Closest Av.'' refers to the closest-available-ambulance heuristic, ``Andersson'' refers to the policy in \citet{ande:07}, ``Lee'' refers to the first modified policy in \cite{lees:11}, ``Mayorga'' refers to the policy in \citet{mayo:13} that applies the heuristic of \cite{band:12} if there is an ambulance available in the same district as the emergency, and that dispatches an ambulance from another district using a preference list of ambulances if there is no ambulance available in the same district as the emergency, ``Bandara'' refers to the policy in \citet{band:14}, and ``Jagtenberg'' refers to the DMEXCLP policy in \citet{jagt:17a,jagt:17b}.
The columns headed ``NR'' refer to the policies as described in the respective papers.
In addition to these policies, we also implemented roll-out policies that use these heuristics in the following way.
Similar to the arc-based and itinerary-based policies, the roll-out policies are based on a two-stage formulation with $15$~second-stage scenarios (independent from the simulated scenario) with a time horizon length of $T = 2$~hours for each second-stage scenario.
The first-stage problem of the roll-out policies are exactly the same as the first-stage problem of the arc-based and itinerary-based policies.
For the second-stage problems, instead of solving the second-stage optimization problems of the arc-based and itinerary-based policies, the roll-out policies use the heuristics (closest-available-ambulance heuristic, \ldots, Jagtenberg DMEXCLP heuristic) to make second-stage decisions, and these second-stage objective values are used to solve the first-stage problems in the same way as for the arc-based and itinerary-based policies.
The columns headed ``RO'' refer to these roll-out policies.
In these columns, each row contains two subrows, with the first subrow (``RT'' for actual response time and ``PRT'' for penalized response time) giving the average objective value, and the second subrow (``SDIFF'') giving the number of standard deviations of the difference between the average objective value of the reference policy and the average objective value of the policy in the particular column (negative values means that the reference policy performed better on average).
The objective values reported in Tables~\ref{tbl:enum_vs_all_actual_random_nro} and~\ref{tbl:arc_vs_all_actual_random_nro} are the average actual response times, that is, the time elapsed from the moment an emergency call is received until an ambulance arrives at the emergency scene, which is the same performance metric considered in many papers on ambulance dispatching.
The objective values reported in Tables~\ref{tbl:enum_vs_all_random_nro} and~\ref{tbl:arc_vs_all_random_nro} are the average penalized response times, which take into account the response time, the type of emergency, as well as the quality of care resulting from the type of ambulance that services the type of emergency.
The penalized response time for each emergency is computed as follows: let $w$ denote the response time (in seconds) of an emergency, and as before let $c$ denote the type of an emergency.
Response time is not equally important for all emergency types, and therefore the response time for an emergency of type~$c$ is penalized with a response time coefficient $\theta_{c}$.
In addition, different emergency types need ambulances and crews with different capabilities, and therefore the objective includes a quality of care coefficient $M_{ac}$ if an ambulance of type~$a$ is used to serve a type~$c$ emergency.
The penalized response time is given by $P(a,c,w) = \theta_{c} w + 3600 M_{ac}$.
For the numerical results, we used $\theta_{c} = 1$ for $c \in \{\mathrm{A,C}\}$ (low time urgency emergencies), and $\theta_{c} = 4$ for $c \in \{\mathrm{B,D}\}$ (high time urgency emergencies).
We used the quality of care coefficients $M_{ac}$ given in Table~\ref{tbl:compatibility_matrix}.

\begin{table}[]
\centering
\begin{tabular}{|c|c|c|c|c|}
\hline
   & A: Low, no pref. & B: High, no pref. & C: Low, ALS pref. & D: High, ALS pref. \\ \hline
ALS & 1 & 1 & 0 & 0 \\ \hline
BLS & 0 & 0 & 4 & 4 \\ \hline
\end{tabular}
\caption{Quality of care coefficients $M_{ac}$. The columns denote each type of emergency, with its respective time urgency and ambulance type preference, and the rows represent the type of ambulance.}
\label{tbl:compatibility_matrix}
\end{table}




\begin{table}[]
\centering
\scalebox{0.7}{
\begin{tabular}{cccc|cl|cl|cl|cl|cl|cl}
\hline
A                    & \begin{tabular}[c]{@{}c@{}}Itinerary \\ PRT\end{tabular} &                      & Arc                   & \multicolumn{2}{c|}{Closest Av.}                    & \multicolumn{2}{c|}{Andersson}  & \multicolumn{2}{c|}{Lee}        & \multicolumn{2}{c|}{Mayorga}    & \multicolumn{2}{c|}{Bandara}    & \multicolumn{2}{c}{Jagtenberg}  \\ \hline
\multicolumn{1}{l}{} & \multicolumn{1}{l}{}                                    & \multicolumn{1}{l}{} & \multicolumn{1}{l|}{} & \multicolumn{1}{l}{RO:} & NR:                       & \multicolumn{1}{l}{RO:} & NR:   & \multicolumn{1}{l}{RO:} & NR:   & \multicolumn{1}{l}{RO:} & NR:   & \multicolumn{1}{l}{RO:} & NR:   & \multicolumn{1}{l}{RO:} & NR:   \\ \cline{5-16} 
\multirow{2}{*}{10}  & \multirow{2}{*}{17368}                                  & PRT:                 & 17525                 & 23844                   & 27220                     & 25201                   & 29262 & 18424                   & 22638 & 45116                   & 47343 & 42057                   & 45897 & 26455                   & 31286 \\
                     &                                                         & SDIFF:               & -0.1                  & -3.4                    & -4.9                      & -4.2                    & -5.3  & -0.6                    & -2.8  & -7.1                    & -7.4  & -6.7                    & -7.2  & -4.8                    & -5.9  \\ \hline
\multirow{2}{*}{12}  & \multirow{2}{*}{14327}                                  & PRT:                 & 15925                 & 21222                   & 26789                     & 21714                   & 24399 & 17936                   & 19706 & 33429                   & 36828 & 38719                   & 39935 & 22138                   & 25743 \\
                     &                                                         & SDIFF:               & -1.5                  & -4.6                    & -4.7                      & -4.8                    & -5.3  & -2.4                    & -4.0  & -6.8                    & -7.1  & -6.8                    & -7.0  & -5                      & -6.2  \\ \hline
\multirow{2}{*}{14}  & \multirow{2}{*}{11076}                                  & PRT:                 & 13204                 & 18604                   & 20625                     & 19783                   & 20982 & 16532                   & 17729 & 28064                   & 29204 & 30579                   & 31450 & 21787                   & 24378 \\
                     &                                                         & SDIFF:               & -2.1                  & -5.4                    & -4.5                      & -5.2                    & -5.6  & -3.7                    & -4.1  & -6.8                    & -6.9  & -7.6                    & -7.9  & -6.8                    & -7.3  \\ \hline
\multirow{2}{*}{16}  & \multirow{2}{*}{10752}                                  & PRT:                 & 14405                 & 18975                   & 20119                     & 17804                   & 19598 & 13128                   & 15386 & 22028                   & 24109 & 24590                   & 26358 & 17996                   & 19822 \\
                     &                                                         & SDIFF:               & -3                    & -4.5                    & -4.1                      & -4.6                    & -4.9  & -2.2                    & -3.0  & -6.2                    & -6.5  & -8                      & -8.2  & -4.6                    & -4.9  \\ \hline
\multirow{2}{*}{18}  & \multirow{2}{*}{11069}                                  & PRT:                 & 12741                 & 15853                   & 16932                     & 16492                   & 18453 & 12981                   & 14167 & 20818                   & 21745 & 22785                   & 24253 & 16201                   & 18564 \\
                     &                                                         & SDIFF:               & -1.5                  & -3.5                    & -3.9                      & -3.6                    & -4.3  & -1.7                    & -3.1  & -5.9                    & -6.1  & -6.3                    & -6.8  & -3.1                    & -4.3  \\ \hline
\multirow{2}{*}{20}  & \multirow{2}{*}{10742}                                  & PRT:                 & 11652                 & 14148                   & 15637                     & 14994                   & 16214 & 13274                   & 15152 & 19126                   & 20855 & 20169                   & 22287 & 15387                   & 16973 \\
                     &                                                         & SDIFF:               & -0.8                  & -2.7                    & -3.2                      & -3.7                    & -4.4  & -1.8                    & -2.2  & -6.1                    & -6.5  & -6.2                    & -6.8  & -3.5                    & -4.1  \\ \hline
\multirow{2}{*}{22}  & \multirow{2}{*}{9059}                                   & PRT:                 & 10819                 & 14034                   & 15250                     & 14607                   & 16046 & 10870                   & 12348 & 17554                   & 19473 & 18489                   & 20310 & 14675                   & 16508 \\
                     &                                                         & SDIFF:               & -1.7                  & -3.4                    & -3.8                      & -4.2                    & 4.8   & -1.6                    & 2.6   & -4.7                    & -5.1  & -6.3                    & -6.9  & -4.8                    & -5.4  \\ \hline
\multirow{2}{*}{24}  & \multirow{2}{*}{6780}                                   & PRT:                 & 8834                  & 10205                   & 11175                     & 11938                   & 13120 & 11940                   & 12723 & 14353                   & 15077 & 15528                   & 16846 & 15346                   & 16176 \\
                     &                                                         & SDIFF:               & -2.4                  & -3                      & -4.2                      & -4.9                    & -5.2  & -4.3                    & -4.7  & -5.9                    & -6.3  & -6.1                    & -7.0  & -5.7                    & -6.5  \\ \hline
\multirow{2}{*}{26}  & \multirow{2}{*}{5464}                                   & PRT:                 & 6901                  & 14034                   & 11583                     & 11818                   & 12741 & 10505                   & 11447 & 13865                   & 14801 & 14090                   & 14954 & 12919                   & 14252 \\
                     &                                                         & SDIFF:               & -2.5                  & -3.4                    & -3.6                      & -4.7                    & -5.1  & -5.3                    & -5.9  & -7.3                    & -7.6  & -6.8                    & -7.2  & -5.8                    & -6.4  \\ \hline
\multirow{2}{*}{28}  & \multirow{2}{*}{6166}                                   & PRT:                 & 6757                  & 9975                    & 10606                     & 11628                   & 12445 & 9644                    & 10340 & 12923                   & 13583 & 13250                   & 13805 & 12710                   & 13455 \\
                     &                                                         & SDIFF:               & -1.1                  & -2.7                    & -3.5                      & -6                      & -6.5  & -3.4                    & -3.7  & -5.8                    & -6.4  & -6.1                    & -7.0  & -6                      & -6.8  \\ \hline
\multirow{2}{*}{30}  & \multirow{2}{*}{5385}                                   & PRT:                 & 6865                  & 9600                    & 10137                     & 11336                   & 11928 & 10401                   & 11049 & 12610                   & 13216 & 12173                   & 12613 & 12096                   & 13034 \\
                     &                                                         & SDIFF:               & -2.3                  & -4                      & \multicolumn{1}{c|}{-4.4} & -5.6                    & -6.2  & -4.7                    & -5.1  & -6.1                    & -6.4  & -6.3                    & -6.6  & -5.9                    & -6.7  \\ \hline
\end{tabular}}
\caption{Comparison of penalized response times (PRT) between the itinerary-based model and other policies.}
\label{tbl:enum_vs_all_random_nro}
\end{table}
\begin{table}[]
\centering
\scalebox{0.7}{
\begin{tabular}{cccc|cc|cc|cc|cc|cc|cc}
\hline
A                   & \begin{tabular}[c]{@{}c@{}}Itinerary \\ RT\end{tabular} &      & Arc  & \multicolumn{2}{c|}{Closest Av.} & \multicolumn{2}{c|}{Andersson} & \multicolumn{2}{c|}{Lee} & \multicolumn{2}{c|}{Mayorga} & \multicolumn{2}{c|}{Bandara} & \multicolumn{2}{c}{Jagtenberg} \\ \hline
                    &                                                         &      &      & RO:            & NR:             & RO:            & NR:           & RO:         & NR:        & RO:           & NR:          & RO:           & NR:          & RO:            & NR:           \\ \cline{5-16} 
\multirow{2}{*}{10} & \multirow{2}{*}{8532}                                   & RT:  & 7929 & 9482           & 10215           & 10204          & 10905         & 7736        & 8999       & 16799         & 17488        & 15205         & 16920        & 10648          & 11909         \\
                    &                                                         & SDIFF: & 1    & -1.6           & -2.2            & -2.8           & -3.1          & 1.2         & -2.1       & -7.5          & -7.9         & -6.5          & -7.0         & -3.2           & -3.7          \\ \hline
\multirow{2}{*}{12} & \multirow{2}{*}{7131}                                   & RT:  & 6454 & 7876           & 9150            & 7948           & 8995          & 7311        & 7644       & 12597         & 12988        & 13730         & 14857        & 8035           & 8310          \\
                    &                                                         & SDIFF: & 1.4  & -1.3           & -2.6            & -1.6           & -3.1          & -0.3        & -1.6       & -6.3          & -6.8         & -6.6          & -6.8         & -1.7           & -3.0          \\ \hline
\multirow{2}{*}{14} & \multirow{2}{*}{5532}                                   & RT:  & 5534 & 6619           & 6633            & 6791           & 7540          & 6027        & 6731       & 10482         & 10916        & 10489         & 11020        & 7639           & 7938          \\
                    &                                                         & SDIFF: & 0    & -2.1           & -2.4            & -2.3           & -2.9          & -1.1        & -2.4       & -6.4          & -6.8         & -6.4          & -6.7         & -4.1           & -4.4          \\ \hline
\multirow{2}{*}{16} & \multirow{2}{*}{5182}                                   & RT:  & 5838 & 6289           & 6580            & 5599           & 5937          & 4908        & 5478       & 7646          & 8284         & 8300          & 8662         & 5677           & 6551          \\
                    &                                                         & SDIFF: & -1.6 & -2             & -2.2            & -0.9           & -1.6          & 0.7         & -0.3       & -4.4          & -5.5         & -5.8          & -6.7         & -1             & -2.2          \\ \hline
\multirow{2}{*}{18} & \multirow{2}{*}{4963}                                   & RT:  & 4975 & 5464           & 6005            & 5459           & 5777          & 4436        & 4945       & 7235          & 7320         & 7747          & 8008         & 5371           & 5887          \\
                    &                                                         & SDIFF: & 0    & -1.2           & -2.3            & -1.1           & -2.3          & 1.4         & 0.1        & -4.7          & -5.4         & -4.9          & -5.9         & -0.9           & -1.7          \\ \hline
\multirow{2}{*}{20} & \multirow{2}{*}{4483}                                   & RT:  & 4047 & 4807           & 5365            & 4425           & 5815          & 4683        & 5235       & 6294          & 6524         & 6473          & 6662         & 4562           & 4816          \\
                    &                                                         & SDIFF: & 1.1  & -0.8           & -1.8            & 0.2            & -2.4          & -0.4        & -2.3       & -4.2          & -5.2         & -5            & -5.1         & -0.2           & -0.8          \\ \hline
\multirow{2}{*}{22} & \multirow{2}{*}{4564}                                   & RT:  & 3785 & 4407           & 5486            & 4392           & 4818          & 3442        & 4550       & 5788          & 6183         & 5940          & 6464         & 4154           & 4630          \\
                    &                                                         & SDIFF: & 1.5  & 0.3            & -1.6            & 0.4            & -0.4          & 2.4         & 0.1        & -2.3          & -2.7         & -2.3          & -2.5         & 1              & -0.6          \\ \hline
\multirow{2}{*}{24} & \multirow{2}{*}{3356}                                   & RT:  & 3228 & 3131           & 5064            & 3386           & 3847          & 3612        & 4038       & 4447          & 5015         & 4568          & 4684         & 4439           & 4986          \\
                    &                                                         & SDIFF: & 0.3  & 0.6            & -2.4            & -0.1           & -0.7          & -0.6        & -0.7       & -2.6          & -3.0         & -2.8          & -3.2         & -2.1           & -3.1          \\ \hline
\multirow{2}{*}{26} & \multirow{2}{*}{2978}                                   & RT:  & 3308 & 3174           & 4667            & 3453           & 3916          & 2958        & 3301       & 4456          & 4989         & 4348          & 4622         & 3758           & 4187          \\
                    &                                                         & SDIFF: & -0.8 & -0.4           & -2.0            & -1.1           & -2.5          & 0.1         & -0.4       & -3.6          & -3.8         & -3.2          & -3.4         & -2.0           & -3.0          \\ \hline
\multirow{2}{*}{28} & \multirow{2}{*}{2784}                                   & RT:  & 2870 & 2492           & 4465            & 2967           & 3117          & 2757        & 3149       & 4141          & 4568         & 3959          & 3635         & 3350           & 3732          \\
                    &                                                         & SDIFF: & -0.2 & 0.8            & -2.6            & -0.6           & -2.3          & 0.1         & -1.8       & -3            & -3.3         & -3.4          & -2.9         & -1.9           & -2.5          \\ \hline
\multirow{2}{*}{30} & \multirow{2}{*}{2833}                                   & RT:  & 2889 & 2499           & 4208            & 3149           & 3667          & 2725        & 2887       & 3513          & 3633         & 3639          & 3111         & 3083           & 3385          \\
                    &                                                         & SDIFF: & -0.2 & 0.8            & -2.1            & -0.7           & -1.6          & 0.3         & -0.1       & -1.5          & -3.1         & -1.9          & -1.5         & -0.6           & -1.8          \\ \hline
\end{tabular}}
\caption{Comparison of actual response times (in seconds) between the itinerary-based model and other policies.}
\label{tbl:enum_vs_all_actual_random_nro}
\end{table}

\begin{table}[]
\centering
\scalebox{0.7}{
\begin{tabular}{cccc|cl|cl|cl|cl|cl|cl}
\hline
A                    & \begin{tabular}[c]{@{}c@{}}Arc\\ PRT\end{tabular} &                      & Itinerary             & \multicolumn{2}{c|}{Closest Av.}                    & \multicolumn{2}{c|}{Andersson} & \multicolumn{2}{c|}{Lee}       & \multicolumn{2}{c|}{Mayorga}   & \multicolumn{2}{c|}{Bandara}   & \multicolumn{2}{c}{Jagtenberg} \\ \hline
\multicolumn{1}{l}{} &                                                  & \multicolumn{1}{l}{} & \multicolumn{1}{l|}{} & \multicolumn{1}{l}{RO:} & NR:                       & \multicolumn{1}{l}{RO:} & NR:  & \multicolumn{1}{l}{RO:} & NR:  & \multicolumn{1}{l}{RO:} & NR:  & \multicolumn{1}{l}{RO:} & NR:  & \multicolumn{1}{l}{RO:} & NR:  \\ \cline{5-16} 
\multirow{2}{*}{10}  & \multirow{2}{*}{17525}                           & PRT:                 & 17368                 & 23844                   & 27220                 & 25201                   & 29262 & 18424                   & 22638 & 45116                   & 47343 & 42057                   & 45897 & 26455                   & 31286 \\
                     &                                                  & SDIFF:                 & 0.1                   & -4.2                    & -4.9                      & -4.5                    & -4.8  & -0.6                    & -2.2 & -6.9                    & -7.4 & -6.6                    & -7.1 & -4.6                    & -5.2 \\ \hline
\multirow{2}{*}{12}  & \multirow{2}{*}{15925}                           & PRT:                 & 14327                 & 21222                   & 26789                 & 21714                   & 24399 & 17936                   & 19706 & 33429                   & 36828 & 38719                   & 39935 & 22138                   & 25743 \\
                     &                                                  & SDIFF:                 & 1.5                   & -3.5                    & -4.9                      & -4.2                    & -5.0 & -1.4                    & -3.0 & -6.1                    & -6.6 & -6.3                    & -6.6 & -3.5                    & -4.2 \\ \hline
\multirow{2}{*}{14}  & \multirow{2}{*}{13204}                           & PRT:                 & 11076                 & 18604                   & 20625                 & 19783                   & 20982 & 16532                   & 17729 & 28064                   & 29204 & 30579                   & 31450 & 21787                   & 24378 \\
                     &                                                  & SDIFF:                 & 2.1                   & -4.2                    & -4.8                      & -4.3                    & -5.0 & -2.3                    & -2.9 & -6.4                    & -6.8 & -7.4                    & -7.7 & -5.7                    & -7.0 \\ \hline
\multirow{2}{*}{16}  & \multirow{2}{*}{14405}                           & PRT:                 & 10752                 & 18975                   & 20119                 & 17804                   & 19598 & 13128                   & 15386 & 22028                   & 24109 & 24590                   & 26358 & 17996                   & 19822 \\
                     &                                                  & SDIFF:                 & 3                     & -2.3                    & -3.1                      & -2                      & -3.3 & 1                       & -1.8 & -4                      & -5.8 & -5.6                    & -6.6 & -2.1                    & -2.9 \\ \hline
\multirow{2}{*}{18}  & \multirow{2}{*}{12741}                           & PRT:                 & 11069                 & 15853                   & 16932                 & 16492                   & 18453 & 12981                   & 14167 & 20818                   & 21745 & 22785                   & 24253 & 16201                   & 18564 \\
                     &                                                  & SDIFF:                 & 1.5                   & -2.3                    & -3.9                      & -2.7                    & -3.8 & -0.2                    & -1.4 & -4.7                    & -5.0 & -5.7                    & -6.4 & -2.2                    & -3.3 \\ \hline
\multirow{2}{*}{20}  & \multirow{2}{*}{11652}                           & PRT:                 & 10742                 & 14148                   & 15637                 & 14994                   & 16214 & 13274                   & 15152 & 19126                   & 20855 & 20169                   & 22287 & 15387                   & 16973 \\
                     &                                                  & SDIFF:                 & 0.8                   & -1.8                    & -2.9                      & -2.3                    & -3.4 & -1.3                    & -3.0 & -5                      & -5.7 & -4.9                    & -6.8 & -2.8                    & -4.1 \\ \hline
\multirow{2}{*}{22}  & \multirow{2}{*}{10819}                           & PRT:                 & 9059                  & 14034                   & 15250                 & 14607                   & 16046 & 10870                   & 12348 & 17554                   & 19473 & 18489                   & 20310 & 14675                   & 16508 \\
                     &                                                  & SDIFF:                 & 1.7                   & -2.2                    & -3.8                      & -2.6                    & -3.8 & 0                       & -1.8 & -3.8                    & -4.3 & -5.5                    & -6.0 & -2.8                    & -3.9 \\ \hline
\multirow{2}{*}{24}  & \multirow{2}{*}{8834}                            & PRT:                 & 6780                  & 10205                   & 11175                 & 11938                   & 13120 & 11940                   & 12723 & 14353                   & 15077 & 15528                   & 16846 & 15346                   & 16176 \\
                     &                                                  & SDIFF:                 & 2.4                   & -1.3                    & -2.4                      & -2.8                    & -3.4 & -2.8                    & -3.7 & -4.4                    & -4.9 & -4.8                    & -5.9 & -4.3                    & -5.1 \\ \hline
\multirow{2}{*}{26}  & \multirow{2}{*}{6901}                            & PRT:                 & 5464                  & 14034                   & 11583                 & 11818                   & 12741 & 10505                   & 11447 & 13865                   & 14801 & 14090                   & 14954 & 12919                   & 14252 \\
                     &                                                  & SDIFF:                 & 2.5                   & -2.2                    & -3.4                      & -2.4                    & -2.9 & -2.1                    & -3.3 & -4.3                    & -5.4 & -4.0                    & -4.4 & -3.2                    & -4.9 \\ \hline
\multirow{2}{*}{28}  & \multirow{2}{*}{6757}                            & PRT:                 & 6166                  & 9975                    & 10606                 & 11628                   & 12445 & 9644                    & 10340 & 12923                   & 13583 & 13250                   & 13805 & 12710                   & 13455 \\
                     &                                                  & SDIFF:                 & 1.1                   & -2.0                    & -3.2                      & -3.7                    & -4.3 & -2.4                    & -2.9 & -5.1                    & -5.8 & -5.4                    & -5.7 & -4.9                    & -5.3 \\ \hline
\multirow{2}{*}{30}  & \multirow{2}{*}{6865}                            & PRT:                 & 5385                  & 9600                    & 10137                 & 11336                   & 11928 & 10401                   & 11049 & 12610                   & 13216 & 12173                   & 12613 & 12096                   & 13034 \\
                     &                                                  & SDIFF:                 & 2.3                   & -1.7                    & \multicolumn{1}{c|}{-2.8} & -3.1                    & -3.3 & -3.0                    & -3.7 & -4.3                    & -4.4 & -4.1                    & -4.6 & -3.4                    & -4.8 \\ \hline
\end{tabular}}
\caption{Comparison of penalized response times between the arc-based model and other policies.}
\label{tbl:arc_vs_all_random_nro}
\end{table}
\begin{table}[]
\centering
\scalebox{0.7}{
\begin{tabular}{cccc|cc|cc|cc|cc|cc|cc}
\hline
A                   & \begin{tabular}[c]{@{}c@{}}Arc\\ RT\end{tabular} &      & Itinerary & \multicolumn{2}{c|}{Closest Av.} & \multicolumn{2}{c|}{Andersson} & \multicolumn{2}{c|}{Lee} & \multicolumn{2}{c|}{Mayorga} & \multicolumn{2}{c|}{Bandara} & \multicolumn{2}{c}{Jagtenberg} \\ \hline
                    &                                                  &      &           & RO:            & NR:             & RO:            & NR:           & RO:         & NR:        & RO:           & NR:          & RO:           & NR:          & RO:            & NR:           \\ \cline{5-16} 
\multirow{2}{*}{10} & \multirow{2}{*}{7929}                            & RT:  & 8532      & 9482           & 10215           & 10204          & 10905         & 7736        & 8999       & 16799         & 17488        & 15205         & 16920        & 10648          & 11909         \\
                    &                                                  & SDIFF: & -1        & -2.3           & -2.7            & -3.4           & -3.7          & 0.4         & -1.6       & -7.5          & -7.8         & -6.8          & -7.1         & -4.0           & -4.4          \\ \hline
\multirow{2}{*}{12} & \multirow{2}{*}{6454}                            & RT:  & 7131      & 7876           & 9150            & 7948           & 8995          & 7311        & 7644       & 12597         & 12988        & 13730         & 14857        & 8035           & 8310          \\
                    &                                                  & SDIFF: & -1.4      & -3.0           & -4.7            & -3.3           & -4.6          & -2          & -2.9       & -6.5          & -6.9         & -6.8          & -7.4         & -2.7           & -3.3          \\ \hline
\multirow{2}{*}{14} & \multirow{2}{*}{5534}                            & RT:  & 5532      & 6619           & 6633            & 6791           & 7540          & 6027        & 6731       & 10482         & 10916        & 10489         & 11020        & 7639           & 7938          \\
                    &                                                  & SDIFF: & 0         & -2.4           & -3.4            & -2.6           & -4.0          & -1          & -1.9       & -6.1          & -7.3         & -6.5          & -7.4         & -4.1           & -4.6          \\ \hline
\multirow{2}{*}{16} & \multirow{2}{*}{5838}                            & RT:  & 5182      & 6289           & 6580            & 5599           & 5937          & 4908        & 5478       & 7646          & 8284         & 8300          & 8662         & 5677           & 6551          \\
                    &                                                  & SDIFF: & 1.6       & -0.8           & -1.6            & 0.9            & -1.3          & 2.4         & 0.2        & -3.0          & -3.6         & -4.0          & -4.3         & 0.3            & -1.2          \\ \hline
\multirow{2}{*}{18} & \multirow{2}{*}{4975}                            & RT:  & 4963      & 5464           & 6005            & 5459           & 5777          & 4436        & 4945       & 7235          & 7320         & 7747          & 8008         & 5371           & 5887          \\
                    &                                                  & SDIFF: & 0         & -1.1           & -2.3            & -1             & -2.5          & 1.2         & -0.1       & -4.1          & -4.6         & -4.7          & -5.7         & -0.7           & -2.0          \\ \hline
\multirow{2}{*}{20} & \multirow{2}{*}{4047}                            & RT:  & 4483      & 4807           & 5365            & 4425           & 5815          & 4683        & 5235       & 6294          & 6524         & 6473          & 6662         & 4562           & 4816          \\
                    &                                                  & SDIFF: & -1.1      & -2.4           & -3.3            & -1.1           & -2.7          & -1.9        & -3.4       & -4.6          & -5.1         & -5.4          & -5.9         & -1.5           & -2.7          \\ \hline
\multirow{2}{*}{22} & \multirow{2}{*}{3785}                            & RT:  & 4564      & 4407           & 5486            & 4392           & 4818          & 3442        & 4550       & 5788          & 6183         & 5940          & 6464         & 4154           & 4630          \\
                    &                                                  & SDIFF: & -1.5      & -1.3           & -2.0            & -1.5           & -2.3          & 0.9         & -2.5       & -3.7          & -4.6         & -4.0          & -4.2         & -0.9           & -2.5          \\ \hline
\multirow{2}{*}{24} & \multirow{2}{*}{3228}                            & RT:  & 3356      & 3131           & 5064            & 3386           & 3847          & 3612        & 4038       & 4447          & 5015         & 4568          & 4684         & 4439           & 4986          \\
                    &                                                  & SDIFF: & -0.3      & 0.3            & -1.7            & -0.5           & -2.3          & -1.1        & -2.5       & -4.3          & -5.1         & -4.5          & -4.6         & -2.9           & -3.5          \\ \hline
\multirow{2}{*}{26} & \multirow{2}{*}{3308}                            & RT:  & 2978      & 3174           & 4667            & 3453           & 3916          & 2958        & 3301       & 4456          & 4989         & 4348          & 4622         & 3758           & 4187          \\
                    &                                                  & SDIFF: & 0.8       & 0.3            & -1.2            & -0.4           & -2.5          & 1           & -0.1       & -2.9          & -3.4         & -2.3          & -3.2         & -1             & -2.3          \\ \hline
\multirow{2}{*}{28} & \multirow{2}{*}{2870}                            & RT:  & 2784      & 2492           & 4465            & 2967           & 3117          & 2757        & 3149       & 4141          & 4568         & 3959          & 3635         & 3350           & 3732          \\
                    &                                                  & SDIFF: & 0.2       & 0.9            & -1.7            & -0.3           & -2.5          & 0.4         & -0.8       & -3.4          & -4.4         & -3.2          & -2.8         & -1.5           & -2.9          \\ \hline
\multirow{2}{*}{30} & \multirow{2}{*}{2889}                            & RT:  & 2833      & 2499           & 4208            & 3149           & 3667          & 2725        & 2887       & 3513          & 3633         & 3639          & 3111         & 3083           & 3385          \\
                    &                                                  & SDIFF: & 0.2       & 1.1            & -1.4            & -0.7           & -2.2          & 0.5         & 0.1        & -1.8          & -2.2         & -2.5          & -3.3         & -0.6           & -2.5          \\ \hline
\end{tabular}}
\caption{Comparison of actual response times (in seconds) between the arc-based model and other policies.}
\label{tbl:arc_vs_all_actual_random_nro}
\end{table}

The results show that both the arc-based policy and the itinerary-based policy performed much better than the basic versions of the other policies (in the ``NR'' columns) in terms of the penalized response time performance metric, and that most of the differences were statistically significant.
The arc-based and itinerary-based policies also performed better than the roll-out policies (in the ``RO'' columns) in terms of the penalized response time performance metric, but the differences were smaller than for the basic versions of the other policies, because the roll-out policies uniformly outperformed the basic versions of the other policies.
Thus, part of the improvement of the arc-based policy and the itinerary-based policy over the other policies was due to the look-ahead feature of the arc-based and itinerary-based policies that was shared with the roll-out policies, and part was due to better modeling of the future costs in the second stage.

The arc-based policy and the itinerary-based policy also performed better than the basic versions of the other policies in terms of the actual response time performance metric, but a few of the differences were not statistically significant.
The performance of the arc-based and itinerary-based policies in terms of the actual response time performance metric was similar to that of many of the roll-out policies, with many of the differences not being statistically significant.

\subsection{Effect of Second-Stage Time Horizon on Proposed Policies}

The numerical results reported in Section~\ref{sec:policy comparison} used a time horizon length of $T = 2$~hours for each second-stage scenario.
We also evaluated the effect of other time horizon lengths on the performance metrics.
Tables~\ref{tab:t_variation_pen} and~\ref{tab:t_variation_actual} compare the average penalized response times and the average actual response times respectively for the itinerary-based and arc-based policies when using time horizon lengths of $T = 2$, $3$, and~$4$~hours for the second stage scenarios, for the setting with $20$~ambulance stations.
The tables follow the same format from previous results.
In most cases, the average penalized response times and the average actual response times decreased as the time horizon length increased, but the differences were not significant.

\begin{table}[]
\centering
\begin{tabular}{ccccc}
\hline
H         & 2     &       & 3     & 4    \\ \hline
Itinerary PRT & 11069 & PRT    & 8898  & 9791 \\
          &       & SDIFF & 1.7   & 0.9  \\ \hline
Arc PRT       & 11652 & PRT    & 11230 & 9975 \\
          &       & SDIFF & 0.4   & 1.1  \\ \hline
\end{tabular}
\caption{Penalized response times for itinerary and arc models with different second stage time horizons. The number of ambulance stations was set to 20.}
\label{tab:t_variation_pen}
\end{table}

\begin{table}[]
\centering
\begin{tabular}{ccccc}
\hline
H                    & 2                    &       & 3    & 4    \\ \hline
Itinerary            & 4487                 & RT    & 4211 & 5116 \\
\multicolumn{1}{l}{} & \multicolumn{1}{l}{} & SDIFF & 0.4  & -1.1 \\ \hline
Arc                  & 4047                 & RT    & 3414 & 3396 \\
\multicolumn{1}{l}{} & \multicolumn{1}{l}{} & SDIFF & 0.6  & 0.6  \\ \hline
\end{tabular}
\caption{Actual response times for itinerary and arc models with different second stage time horizons. The number of ambulance stations was set to 20.}
\label{tab:t_variation_actual}
\end{table}

\subsection{Computation Times of Proposed Policies}

\begin{figure}[h]
    \centering
\includegraphics[scale=0.30]{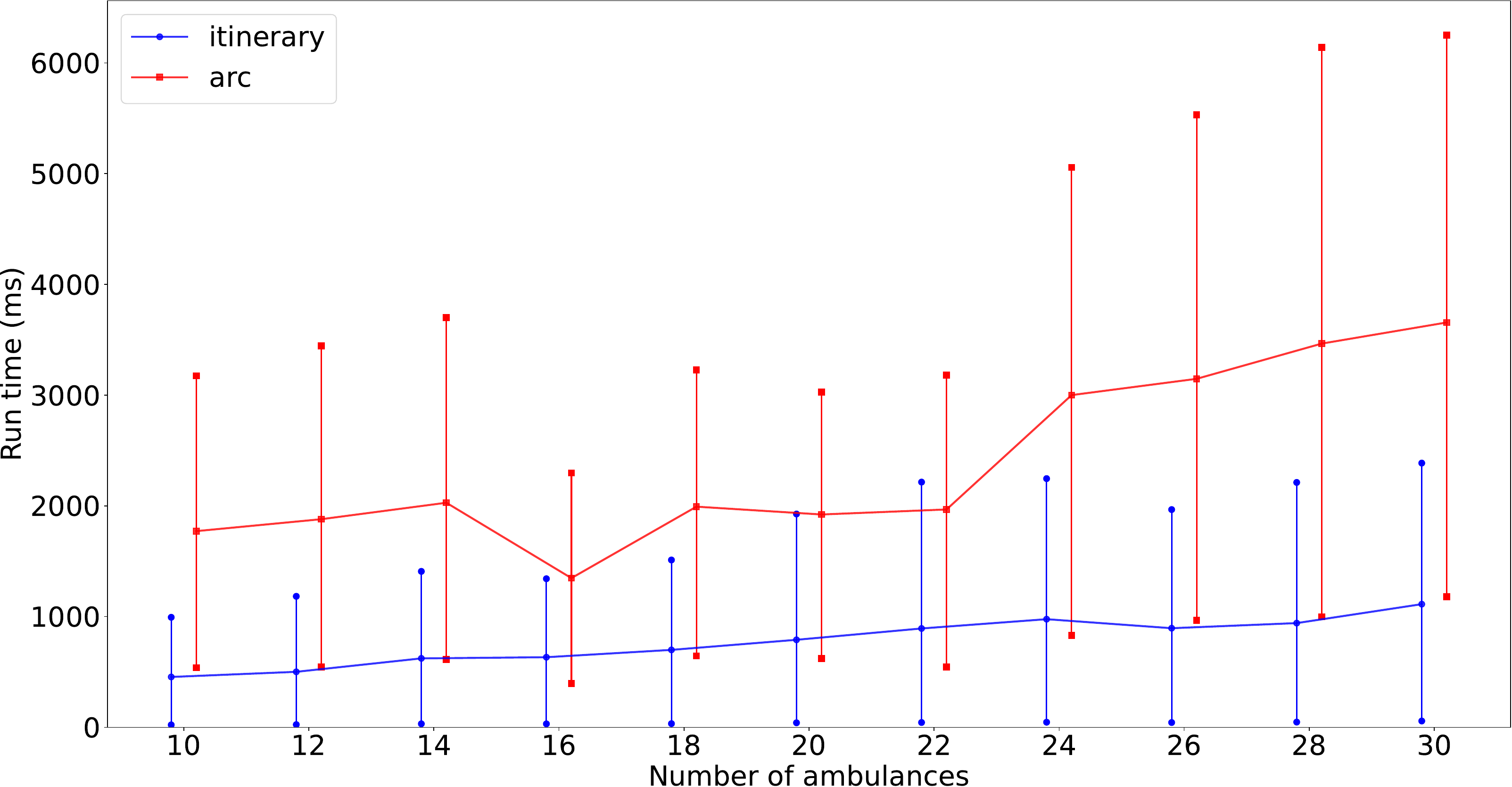}
    \caption{The $0.1$-quantile, the average, and the $0.9$-quantile of computation times for the itinerary-based and arc-based models.}
    \label{figure_runtimes}
\end{figure}

Computation time is an important consideration when using ambulance dispatch policies.
Figure~\ref{figure_runtimes} reports the $0.1$-quantile, the average, and the $0.9$-quantile, of the computation times (in milliseconds) to solve problem~\eqref{eqn:itinerary opt problem objective}--\eqref{eqn:itinerary opt problem each emergency one visit} for the itinerary-based policy and problems~\eqref{eqn:arc second-stage select} and~\eqref{eqn:arc second-stage reassign} for the arc-based policy, as a function of the number of ambulances and ambulance stations.
As can be seen from this figure, with the given settings the computation for both policies take at most a few seconds, and the computations for the itinerary-based policy take less time.

\section{Conclusion}
\label{sec:conc}

The main contribution of this paper is a model of ambulance fleet operations.
The model incorporates many important aspects of ambulance fleet operations, most of which have been ignored in the operations research literature, including the following:
\begin{enumerate}
\item
The model incorporates different emergency types and the consequences of emergency types, such as required ambulance and crew capabilities, hospital requirements, and the marginal value of response time.
\item
The model facilitates different ambulance and crew capabilities.
\item
The model allows hospital choice, taking into account the emergency type and the location of the emergency relative to hospitals.
\item
The model makes provision for a queue of waiting emergencies.
\item
The model incorporates both ambulance selection decisions as well as ambulance reassignment decisions.
\item
The model allows ambulances on their way to a station to be dispatched to an emergency.
\item
The model takes into account the future consequences of dispatch decisions.
\end{enumerate}
We proposed two optimization methods for ambulance dispatch decisions.
The performance of the two methods were compared with six approaches proposed in the literature, using a simulation calibrated with data of the Rio de Janeiro EMS.
The optimization methods performed much better than the other approaches, especially when taking the type of the emergency and the ambulance/crew requirements into account.

Much work remains to be done to improve EMS operations.
The relation between response time and patient outcome is not well understood for most emergency types, with some exceptions such as cardiac arrest and stroke.
A better understanding of the relation between response time and patient outcome for all major emergency types is of great importance for improving EMS performance.
Also, a better understanding is needed of the impact of treatments on patient outcomes, especially longer term outcomes such as survival to hospital discharge and cerebral performance versus shorter term outcomes such as survival to hospital admission \citep[see, for example,][]{ashb:22,mene:23}.
In addition, there is great variability among EMSs in terms of ambulance capabilities and personnel qualifications, varying from EMSs with all ambulances being ALS units and all responders being paramedics, to EMSs that use a mixture of BLS, ILS, ALS units, and even motorcycles and bicycles, and personnel including EMTs, paramedics, nurses, and physicians.
This suggests that a better understanding is needed of the mixture of ambulance and crew capabilities that provides the best outcomes for a given investment \cite[see, for example,][]{nich:96,pers:03}.


\ACKNOWLEDGMENT{The first author's research  was partially supported by an FGV grant.}


\bibliographystyle{informs2014} 
\bibliography{ambulance} 


%
%
%

\newpage

\begin{APPENDIX}{Electronic Companion}

\section{Comparison of Ambulance Dispatching Literature}

Table~\ref{tab:literature summary} compares the ambulance dispatching policies proposed in the literature with the policies considered in this paper.

\begin{center}
\begin{landscape}
\begin{longtable}{|P{0.16\linewidth}|P{0.14\linewidth}|c|P{0.28\linewidth}|P{0.2\linewidth}|}
\hline
 & & Different & & \\
Paper & Emergency Types & Ambulance/Crew & Ambulance Selection Policies & Reassignment Policies \\
 & Considered & Capabilities? & & \\ \hline
\citet{ande:07} & Priority 1, 2, 3 & \ding{55} & Priority~1: dispatch~closest available~ambulance; $\rule{30mm}{0mm}$ Priorities 2, 3: dispatch available ambulance sufficiently close that will maximize minimum preparedness metric of remaining available ambulances & Minimize relocation travel time to raise minimum preparedness metric above threshold \\ \hline
\citet{lees:11} & Single type & \ding{55} & Dispatch available ambulance sufficiently close that will maximize minimum composite preparedness metric of remaining available ambulances; \hspace{100mm} Emergencies are queued if no ambulance available & No relocation; $\rule{20mm}{0mm}$ Upon completion of task, ambulance dispatched to closest emergency in queue if any, otherwise stay in place \\ \hline
\citet{schm:12} & Priorities were mentioned, but all requests were assumed to be of the same priority & \ding{55} & Ambulances dispatched according to approximate dynamic programming solution; $\rule{40mm}{0mm}$ Emergencies are queued if no ambulance \mbox{available} & Ambulances sent to ambulance stations according to approximate dynamic programming solution \\ \hline
\citet{band:12} & Priority 1, 2 & \ding{55} & Ambulances dispatched \mbox{according} to continuous-time MDP solution (solved for $2$,$3$~zones and $2$~ambulances); \hspace{100mm} No queueing of emergencies: Either an ambulance is dispatched immediately or the emergency is served by ``outside'' resources & Upon completion of task, ambulance sent to its home station \\ \hline
\citet{mayo:13} & Priority 1, 2 & \ding{55} & Ambulances dispatched according to one of four district-based dispatching heuristics; $\rule{35mm}{0mm}$ No queueing of emergencies: Either an ambulance is dispatched immediately or the emergency is served by ``outside'' resources & Upon completion of task, ambulance sent to its home station \\ \hline
\citet{band:14} & Priority 1, 2 & \ding{55} & Ambulances dispatched according to closest-available-ambulance rule for Priority~1 emergencies, and a contingency table for Priority~2 emergencies; No queueing of emergencies: Either an ambulance is dispatched immediately or the emergency is served by ``outside'' resources & Upon completion of task, ambulance sent to its home station \\ \hline
\citet{lisay:16} & Priority 1, 2, 3 & \ding{55} & A tabu search algorithm is used to choose a home station for each ambulance; $\rule{50mm}{0mm}$ Priority 1: dispatch closest available ambulance; $\rule{35mm}{0mm}$ Priorities 2, 3: dispatch ambulance sufficiently close with least utilization \hspace{100mm} Emergencies are missed if no ambulance available & Upon completion of task, ambulance sent to its home station \\ \hline
\citet{jagt:17a,jagt:17b} & Single type & \ding{55} & Dispatch ambulance sufficiently close that will maximize expected coverage metric of \mbox{remaining} available ambulances; $\rule{45mm}{0mm}$ Emergencies are queued if no ambulance available & Upon completion of task, ambulance dispatched to oldest emergency in queue if any, otherwise ambulance sent to its home station \\ \hline
This paper & General, as specified by the user, such as MPDS or APCO or NEMSIS chief complaint types or emergency codes & \ding{51} & Ambulance dispatch decisions determined by two-stage stochastic program & Ambulance reassignment decisions determined by two-stage stochastic program \\ \hline
\caption{Summary of papers that consider ambulance dispatch decisions other than the closest-available-ambulance rule.}
\label{tab:literature summary}
\end{longtable}
\end{landscape}
\end{center}

\section{Summary of Notation for Model Input Parameters}

The input parameters of the models are summarized in Table~\ref{pbparameters}.

\begin{table}[h]
\centering
\begin{tabular}{|l|c|}
\hline
~Length of the time horizon & T \\
\hline
~Set of locations & $\mathcal{L}$ \\
\hline
~Set of emergency types & $\mathcal{C}$ \\
\hline
~Set of ambulances & $\mathcal{A}$ \\
\hline
~Set of stations & $\mathcal{B}$ \\
\hline
~Set of hospitals & $\mathcal{H}$ \\
\hline
~Set of ambulances that can serve emergency type $c \in \mathcal{C}$ & $\mathcal{A}(c)$ \\
\hline
~Set of candidate hospitals for emergency type $c \in \mathcal{C}$ at location $\ell \in \mathcal{L}$ & $\mathcal{H}(c,\ell)$ \\
\hline
~Set of emergency types that can be served by ambulance $a \in \mathcal{A}$ & $\mathcal{C}(a)$ \\
\hline
~Set of emergencies & $\mathcal{E}$ \\
\hline
~Set of emergencies that can be handled by ambulance $a \in \mathcal{A}$ & $\mathcal{E}(a)$ \\
\hline
~Arrival time of emergency~$e \in \mathcal{E}$ & $t_{1}(e)$ \\
\hline
~Type of emergency~$e \in \mathcal{E}$ & $c(e)$ \\
\hline
~Location of emergency~$e \in \mathcal{E}$ & $\ell(e)$ \\
\hline
~Set of feasible itineraries for ambulance~$a \in \mathcal{A}$ over $[0,T]$ & $\mathcal{I}(a)$ \\
\hline
~Set of all feasible itineraries over $[0,T]$ & $\mathcal{I}$ \\
\hline
~Type of the emergency at $t=0$ & $c_{0}$ \\
\hline
~Location of the emergency at $t=0$ & $\ell_{0}$ \\
\hline
~Ambulance that becomes available at $t=0$ & $a_{0}$ \\
\hline
~Hospital where the ambulance becomes available at $t=0$ & $h_{0}$ \\
\hline
~Number of emergencies of type $c \in \mathcal{C}$ at location $\ell \in \mathcal{L}$ in queue at $t=0$ & $C_{0}(c,\ell)$ \\
\hline
~Indicator whether ambulance $a \in \mathcal{A}$ available at station $b \in \mathcal{B}$ at $t=0$ & $A_{0}(a,b)$ \\
\hline
~Indicator whether ambulance $a \in \mathcal{A}$ available at hospital $h \in \mathcal{H}$ at $t=0$ & $A_{0}(a,h)$ \\
\hline
\begin{tabular}{l}
Location at time $t+1$ of ambulance~$a$ \\
which is at $\ell_{1}$ at time $t$ on its way to $\ell_{2}$
\end{tabular}
& $L(t,a,\ell_{1},\ell_{2})$ \\
\hline
\begin{tabular}{l}
Set of intermediate locations that can be reached at time~$t$ \\
by ambulance~$a$ going from a hospital to station~$b$
\end{tabular}
& $\mathcal{L}(t,a,b)$ \\
\hline
~Indicator whether ambulance~$a$ at location~$\ell_{1}$ at $t=0$ going to station~$b$
& $A_{0}(a,\ell_{1},b)$ \\
\hline
\begin{tabular}{l}
Indicator whether ambulance~$a$ at location~$\ell_{1}$ at $t=0$ going to \\
an emergency of type $c$ at location $\ell$ and then going to hospital~$h$
\end{tabular}
& $A_{0}(c,a,\ell_{1},\ell,h)$ \\
\hline
\begin{tabular}{l}
Indicator whether ambulance~$a$ at location~$\ell_{1}$ at $t=0$ traveling to \\
hospital~$h$ with emergency type~$c$ patient after on-site emergency care \\
has been provided
\end{tabular}
& $A_{0}(c,a,\ell_{1},h)$ \\
\hline
~Number of emergencies of type~$c$ for period~$t$ and location~$\ell$
& $\lambda(t,c,\ell)$ \\
\hline
\begin{tabular}{l}
Time for ambulance~$a$ to go from location~$\ell_{1}$ at time~$t$ to location~$\ell$, \\
provide on-site emergency care for an emergency of type~$c$ at location~$\ell$, \\
travel from location~$\ell$ to hospital~$h$ and deliver the patient to hospital~$h$
\end{tabular}
& $\tau(t,c,a,\ell_{1},\ell,h)$ \\
\hline
\begin{tabular}{l}
Time for ambulance~$a$ to travel with a patient of type~$c$ \\
from location~$\ell_{1}$ at time~$t=0$, after on-site emergency care has been provided, \\
to hospital~$h$ and deliver the patient at hospital~$h$
\end{tabular}
& $\tau_{0}(c,a,\ell_{1},h)$ \\
\hline
~Maximum number of ambulances at station $b$ & ${A}_{\max}(b)$ \\
\hline
\end{tabular}
\caption{Model input parameters}
\label{pbparameters}
\end{table}

\newpage

\section{Column Generation Algorithm}
\label{sec:column generation details}

\noindent
The dual variables for problem~\eqref{eqn:arc second-stage select} are as follows:
\begin{itemize}
\item
For each $t = 1,\ldots,T$, $a \in \mathcal{A}$, $b \in \mathcal{B}$, let $\beta_{t}(a,b)$ denote the dual variable associated with the flow balance constraint~(At1) at station~$b$.
\item
For each $t = 1,\ldots,T$, $a \in \mathcal{A}$, $h \in \mathcal{H}$, let $\alpha_{t}(a,h)$ denote the dual variable associated with the flow balance constraint~(At2) at hospital~$h$.
\item
For each $t = 1,\ldots,T$, $a \in \mathcal{A}$, $b \in \mathcal{B}$, $\ell_{1} \in \mathcal{L}(t,a,b)$, let $\psi_{t}(a,b,\ell_{1})$ denote the dual variable associated with the flow balance constraint~(At3) at location~$\ell_{1}$.
\item
For each $t = 1,\ldots,T$, $c \in \mathcal{C}$, $\ell \in \mathcal{L}$, let $\phi_{t}(c,\ell)$ denote the dual variable associated with the flow balance constraint~(At4) for the queue of emergency type~$c$ at location~$\ell$.
\item
For each $t = 1,\ldots,T+1$, $b \in \mathcal{B}$, let $\nu_{t}(b) \ge 0$ denote the dual variable associated with the capacity constraint~(At5) for station~$b$.
\item
For each $t = 1,\ldots,T$, $a \in \mathcal{A}$, $b \in \mathcal{B}$, $\ell_{1} \in \mathcal{L}(t,a,b)$, let $\theta_{t}(a,b,\ell_{1}) \ge 0$ denote the dual variable associated with the ambulance supply constraint~(At6) for location~$\ell_{1}$.
\item
For each $t = 1,\ldots,T$, $a \in \mathcal{A}$, $b \in \mathcal{B}$, let $\gamma_{t}(a,b) \ge 0$ denote the dual variable associated with the ambulance supply constraint~(At7) for station~$b$.
\item
For each $t = 1,\ldots,T$, $c \in \mathcal{C}$, $\ell \in \mathcal{L}$, let $\zeta_{t}(c,\ell) \ge 0$ denote the dual variable associated with the ambulance dispatch constraint~(At8).
\end{itemize}

\noindent
For problem~\eqref{eqn:arc second-stage select}, the Lagrangian relaxation is to minimize
\begin{eqnarray*}
\label{eqn:Lagrangianobjectivecall1}
& & \displaystyle{\sum_{t=1}^{T} \sum_{c \in \mathcal{C}} \sum_{a \in \mathcal{A}(c)} \sum_{b \in \mathcal{B}} \sum_{\ell \in \mathcal{L}} \sum_{h \in \mathcal{H}(c,\ell)} \big[f_{t}(c,a,b,\ell,h) + \beta_{t}(a,b) - \alpha_{t + \tau(t,c,a,b,\ell,h)}(a,h) + \phi_{t}(c,\ell)} \\
& & \hspace{110mm} {} + \gamma_{t}(a,b) + \zeta_{t}(c,\ell)\big] x_{t}(c,a,b,\ell,h) \\
& & {} + \displaystyle{\sum_{t=1}^{T} \sum_{c \in \mathcal{C}} \sum_{a \in \mathcal{A}(c)} \sum_{\ell \in \mathcal{L}} \sum_{h \in \mathcal{H}(c,\ell)} \sum_{h' \in \mathcal{H}} \big[f_{t}(c,a,h',\ell,h) + \alpha_{t}(a,h') - \alpha_{t + \tau(t,c,a,h',\ell,h)}(a,h) + \phi_{t}(c,\ell)} \\
& & \hspace{120mm} {} + \zeta_{t}(c,\ell)\big] x_{t}(c,a,h',\ell,h) \\
& & {} + \displaystyle{\sum_{t=1}^{T} \sum_{c \in \mathcal{C}} \sum_{a \in \mathcal{A}(c)} \sum_{b \in \mathcal{B}} \sum_{\ell_{1} \in \mathcal{L}(t,a,b)} \sum_{\ell \in \mathcal{L}} \sum_{h \in \mathcal{H}(c,\ell)} \big[f_{t}(c,a,\ell_{1},b,\ell,h) + \beta_{t}(a,b) \mathds{1}_{\{L(t,a,\ell_{1},b) = b\}} - \alpha_{t + \tau(t,c,a,\ell_{1},\ell,h)}(a,h)} \\
& & \hspace{25mm} {} + \psi_{t}(a,b,L(t,a,\ell_{1},b)) \mathds{1}_{\{L(t,a,\ell_{1},b) \in \mathcal{L}(t,a,b)\}} + \phi_{t}(c,\ell) + \theta_{t}(a,b,\ell_{1}) + \zeta_{t}(c,\ell)\big] x_{t}(c,a,\ell_{1},b,\ell,h) \\
& & {} + \displaystyle{\sum_{t=1}^{T} \sum_{a \in \mathcal{A}} \sum_{b \in \mathcal{B}} \sum_{h \in \mathcal{H}} \left[f_{t}(a,h,b) + \alpha_{t}(a,h) - \beta_{t}(a,b) \mathds{1}_{\{L(t,a,h,b) = b\}} - \psi_{t}(a,b,L(t,a,h,b)) \mathds{1}_{\{L(t,a,h,b) \in \mathcal{L}(t,a,b)\}}\right] y_{t}(a,h,b)} \\
& & {} + \displaystyle{\sum_{t=1}^{T} \sum_{c \in \mathcal{C}} \sum_{\ell \in \mathcal{L}} \left[g_{t}(c,\ell) - \phi_{t}(c,\ell) + \phi_{t-1}(c,\ell)\right] C_{t}(c,\ell)
+ \sum_{c \in \mathcal{C}} \sum_{\ell \in \mathcal{L}} \left[g_{T+1}(c,\ell) + \phi_{T}(c,\ell)\right] C_{T+1}(c,\ell)} \\
& & {} + \displaystyle{\sum_{t=1}^{T} \sum_{a \in \mathcal{A}} \sum_{b \in \mathcal{B}} \left[g_{t}(a,b) - \beta_{t}(a,b) + \beta_{t-1}(a,b) + \nu_{t}(b) - \gamma_{t}(a,b)\right] A_{t}(a,b)} \\
& & {} + \displaystyle{\sum_{a \in \mathcal{A}} \sum_{b \in \mathcal{B}} \left[g_{T+1}(a,b) + \beta_{T}(a,b) + \nu_{T+1}(b)\right] A_{T+1}(a,b)} \\
& & {} + \displaystyle{\sum_{t=1}^{T} \sum_{a \in \mathcal{A}} \sum_{b \in \mathcal{B}} \sum_{\ell_{1} \in \mathcal{L}(t,a,b)} \left[g_{t}(a,\ell_{1},b) - \beta_{t}(a,b) \mathds{1}_{\{L(t,a,\ell_{1},b) = b\}} + \psi_{t-1}(a,b,\ell_{1})\right.} \\
& & \hspace{50mm} \left. {} - \psi_{t}(a,b,L(t,a,\ell_{1},b)) \mathds{1}_{\{L(t,a,\ell_{1},b) \in \mathcal{L}(t,a,b)\}} - \theta_{t}(a,b,\ell_{1})\right] A_{t}(a,\ell_{1},b) \\
& & {} + \displaystyle{\sum_{a \in \mathcal{A}} \sum_{b \in \mathcal{B}} \sum_{\ell_{1} \in \mathcal{L}(T+1,a,b)} \left[g_{T+1}(a,\ell_{1},b) + \psi_{T}(a,b,\ell_{1})\right] A_{T+1}(a,\ell_{1},b)} \\
& & {} + \displaystyle{\sum_{t = 1}^{T} \sum_{a \in \mathcal{A}} \sum_{h \in \mathcal{H}} \alpha_{t}(a,h) \left\{- \sum_{c \in \mathcal{C}(a)} \sum_{\ell_{1} \in \mathcal{L}} \sum_{\{\ell \in \mathcal{L} \, : \, \tau(0,c,a,\ell_{1},\ell,h) = t\}} A_{0}(c,a,\ell_{1},\ell,h)\right.} \\
& & \hspace{50mm} \left. {} - \displaystyle{\sum_{c \in \mathcal{C}(a)} \sum_{\{\ell_{1} \in \mathcal{L} \, : \, \tau_{0}(c,a,\ell_{1},h) = t\}} A_{0}(c,a,\ell_{1},h)}\right\} \\
& & {} - \displaystyle{\sum_{t = 1}^{T} \sum_{c \in \mathcal{C}} \sum_{\ell \in \mathcal{L}} \phi_{t}(c,\ell) \lambda(t,c,\ell) - \sum_{t = 1}^{T+1} \sum_{b \in \mathcal{B}} \nu_{t}(b) {A}_{\max}(b)}
- \displaystyle{\sum_{t = 1}^{T} \sum_{c \in \mathcal{C}} \sum_{\ell \in \mathcal{L}} \zeta_{t}(c,\ell) \left[C_{t}(c,\ell) + \lambda(t,c,\ell)\right]}.
\end{eqnarray*}
Below we consider column generation for three types of primal decision variables: $x_{t}(c,a,b,\ell,h)$, $x_{t}(c,a,h',\ell,h)$, and $x_{t}(c,a,\ell_{1},b,\ell,h)$.

\subsection{Column Generation for Variables $x_{t}(c,a,b,\ell,h)$}

Next, we consider the column generation subproblem to find a variable $x_{t}(c,a,b,\ell,h)$ with the smallest objective coefficient in the Lagrangian relaxation, that is,
\begin{eqnarray}
\min & & \big\{f_{t}(c,a,b,\ell,h) + \beta_{t}(a,b) - \alpha_{t + \tau(t,c,a,b,\ell,h)}(a,h) + \phi_{t}(c,\ell) + \gamma_{t}(a,b) + \zeta_{t}(c,\ell) \; : \; \nonumber \\
& & \hspace{10mm} t \in \{1,\ldots,T\}, c \in \mathcal{C}, a \in \mathcal{A}(c), b \in \mathcal{B}, \ell \in \mathcal{L}, h \in \mathcal{H}(c,\ell)\big\}.
\label{eqn:columngenerationsubproblem3}
\end{eqnarray}
Consider the following simplifying assumptions:
\begin{enumerate}
\item[(B1)]
The dependence of the time $\tau(t,c,a,b,\ell,h)$ on the emergency type~$c$, can be ignored.
Thus this time is denoted with $\tau(t,a,b,\ell,h)$.
\item[(B2)]
The dependence of the cost $f_{t}(c,a,b,\ell,h)$ on the emergency type~$c$, can be ignored.
Thus this cost is denoted with $f_{t}(a,b,\ell,h)$.
\item[(A3)]
The set $\mathcal{A}(c)$ of ambulances to choose from does not depend on the emergency type~$c$.
Thus $\mathcal{A}(c) = \mathcal{A}$ for all~$c$.
\item[(A4)]
The set $\mathcal{H}(c,\ell)$ of hospitals to choose from does not depend on the emergency type~$c$.
Thus this set is denoted with $\mathcal{H}(\ell)$.
\end{enumerate}
Then the calculations can be streamlined as follows:
\begin{enumerate}
\item
For each time~$t$ and emergency location $\ell \in \mathcal{L}$, let
\[
\hat{c}_{t}(\ell) \ \ = \ \
c^*_{t}(\ell) \ \ \in \ \ \argmin\left\{\phi_{t}(c,\ell) + \zeta_{t}(c,\ell) \; : \; c \in \mathcal{C}\right\}
\]
denote the critical emergency type at time~$t$ and location~$\ell$.
\item
For each time~$t$, ambulance~$a \in \mathcal{A}$, station~$b \in \mathcal{B}$, and emergency location $\ell \in \mathcal{L}$, let
\[
\hat{h}_{t}(a,b,\ell) \ \ \in \ \ \argmin\left\{f_{t}(a,b,\ell,h) - \alpha_{t + \tau(t,a,b,\ell,h)}(a,h) \; : \; h \in \mathcal{H}(\ell)\right\}
\]
denote the critical hospital at time~$t$ for ambulance~$a$, station~$b$, and emergency location~$\ell$.
In practice, the set $\mathcal{H}(\ell)$ is small, often a singleton, so the computation of $\hat{h}_{t}(a,b,\ell)$ is quick.
\end{enumerate}
Then the column generation subproblem~(\ref{eqn:columngenerationsubproblem3}) reduces to
\begin{eqnarray}
\min & & \Big\{f_{t}(a,b,\ell,\hat{h}_{t}(a,b,\ell)) + \beta_{t}(a,b) - \alpha_{t + \tau(t,a,b,\ell,\hat{h}_{t}(a,b,\ell))}(a,\hat{h}_{t}(a,b,\ell)) \nonumber \\
& & \hspace{5mm} {} + \phi_{t}(\hat{c}_{t}(\ell),\ell) + \gamma_{t}(a,b) + \zeta_{t}(\hat{c}_{t}(\ell),\ell) \; : \; t \in \{1,\ldots,T\}, a \in \mathcal{A}, b \in \mathcal{B}, \ell \in \mathcal{L}\Big\}.
\label{eqn:columngenerationsubproblem4}
\end{eqnarray}
For any emergency location~$\ell$, an ambulance at a station~$b$ that is close to $\ell$ is more attractive than an ambulance at a station~$b$ that is far from $\ell$.
These observations motivate Algorithm~\ref{alg:columngeneration2} to solve problem~\eqref{eqn:columngenerationsubproblem4}.

\begin{algorithm}
\caption{Column generation algorithm for variables $x_{t}(c,a,b,\ell,h)$}
\begin{algorithmic}[1]
\STATE
For each emergency location $\ell \in \mathcal{L}$, construct a list $\mathcal{B}(\ell)$ of stations $b \in \mathcal{B}$ sorted from closest to $\ell$ to furthest from $\ell$;
\STATE
Find an initial feasible solution, and solve the restricted version of the continuous relaxation of
 problem~\eqref{eqn:arc second-stage select} with the decision variables that are nonzero in the initial feasible solution;
\STATE
optimality$\_$verified $\leftarrow$ false;
\WHILE{not optimality$\_$verified}
    \STATE
    optimality$\_$verified $\leftarrow$ true;
    \STATE
    Use the optimal dual variables for the restricted version of the continuous relaxation of problem~\eqref{eqn:arc second-stage select} to compute $\hat{c}_{t}(\ell)$ and $\hat{h}_{t}(a,b,\ell)$;
	\FOR{$t \in \{1,\ldots,T\}$}
	    \FOR{$a \in \mathcal{A}$}
            \FOR{$\ell \in \mathcal{L}$}
                \STATE
                negative$\_$found $\leftarrow$ false;
	            \FOR{$b \in \mathcal{B}(\ell)$ (from closest to furthest) while not negative$\_$found}
                	\IF{objective value of column generation subproblem~\eqref{eqn:columngenerationsubproblem4} is negative}
                        \STATE
	                    negative$\_$found $\leftarrow$ true;
                        \STATE
                        optimality$\_$verified $\leftarrow$ false;
                        \STATE
                        Add variable $x_{t}(\hat{c}_{t}(\ell),a,b,\ell,\hat{h}_{t}(a,b,\ell))$ to the restricted version of the continuous relaxation of problem~\eqref{eqn:arc second-stage select};
                    \ENDIF
                \ENDFOR
            \ENDFOR
        \ENDFOR
    \ENDFOR
    \IF{not optimality$\_$verified}
        \STATE
        Solve the current restricted version of the continuous relaxation of problem~\eqref{eqn:arc second-stage select};
    \ENDIF
\ENDWHILE
\end{algorithmic}
\label{alg:columngeneration2}
\end{algorithm}

\subsection{Column Generation for Variables $x_{t}(c,a,h',\ell,h)$}

Next, we consider the column generation subproblem to find a variable $x_{t}(c,a,h',\ell,h)$ with the smallest objective coefficient in the Lagrangian relaxation, that is,
\begin{eqnarray}
\min & & \big\{f_{t}(c,a,h',\ell,h) + \alpha_{t}(a,h') - \alpha_{t + \tau(t,c,a,h',\ell,h)}(a,h) + \phi_{t}(c,\ell) + \zeta_{t}(c,\ell) \; : \; \nonumber \\
& & \hspace{5mm} t \in \{1,\ldots,T\}, c \in \mathcal{C}, a \in \mathcal{A}(c), h' \in \mathcal{H}, \ell \in \mathcal{L}, h \in \mathcal{H}(c,\ell)\big\}.
\label{eqn:columngenerationsubproblem5}
\end{eqnarray}
Consider the following simplifying assumptions:
\begin{enumerate}
\item[(C1)]
The dependence of the time $\tau(t,c,a,h',\ell,h)$ on the emergency type~$c$, can be ignored.
Thus this time is denoted with $\tau(t,a,h',\ell,h)$.
\item[(C2)]
The dependence of the cost $f_{t}(c,a,h',\ell,h)$ on the emergency type~$c$, can be ignored.
Thus this cost is denoted with $f_{t}(a,h',\ell,h)$.
\item[(A3)]
The set $\mathcal{A}(c)$ of ambulances to choose from does not depend on the emergency type~$c$.
Thus $\mathcal{A}(c) = \mathcal{A}$ for all~$c$.
\item[(A4)]
The set $\mathcal{H}(c,\ell)$ of hospitals to choose from does not depend on the emergency type~$c$.
Thus this set is denoted with $\mathcal{H}(\ell)$.
\end{enumerate}
Then the calculations can be streamlined as follows:
\begin{enumerate}
\item
For each time~$t$ and emergency location $\ell \in \mathcal{L}$, let
\[
\check{c}_{t}(\ell) \ \ = \ \
c^*_{t}(\ell) \ \ \in \ \ \argmin\left\{\phi_{t}(c,\ell) + \zeta_{t}(c,\ell) \; : \; c \in \mathcal{C}\right\}
\]
denote the critical emergency type at time~$t$ and location~$\ell$.
\item
For each time~$t$, ambulance $a \in \mathcal{A}$, hospital~$h' \in \mathcal{H}$, and emergency location $\ell \in \mathcal{L}$, let
\[
\check{h}_{t}(a,h',\ell) \ \ \in \ \ \argmin\left\{f_{t}(a,h',\ell,h) - \alpha_{t + \tau(t,a,h',\ell,h)}(a,h) \; : \; h \in \mathcal{H}(\ell)\right\}
\]
denote the critical hospital at time~$t$ for ambulance~$a$, hospital~$h'$, and emergency location~$\ell$.
\end{enumerate}
Then the column generation subproblem~(\ref{eqn:columngenerationsubproblem5}) reduces to
\begin{eqnarray}
\min & & \Big\{f_{t}(a,h',\ell,\check{h}_{t}(a,h',\ell)) + \alpha_{t}(a,h') - \alpha_{t + \tau(t,a,h',\ell,\check{h}_{t}(a,h',\ell))}(a,\check{h}_{t}(a,h',\ell)) + \phi_{t}(\check{c}_{t}(\ell),\ell) + \zeta_{t}(\check{c}_{t}(\ell),\ell) \; : \; \nonumber \\
& & \hspace{5mm} t \in \{1,\ldots,T\}, a \in \mathcal{A}, h' \in \mathcal{H}, \ell \in \mathcal{L}
\Big\}.
\label{eqn:columngenerationsubproblem6}
\end{eqnarray}
For any emergency location~$\ell$, an ambulance at a hospital~$h'$ that is close to $\ell$ is more attractive than an ambulance at a hospital~$h'$ that is far from $\ell$.
These observations motivate Algorithm~\ref{alg:columngeneration3} to solve problem~\eqref{eqn:columngenerationsubproblem6}.

\begin{algorithm}
\caption{Column generation algorithm for variables $x_{t}(c,a,h',\ell,h)$}
\begin{algorithmic}[1]
\STATE
For each emergency location $\ell \in \mathcal{L}$, construct a list $\mathcal{H}(\ell)$ of hospitals $h' \in \mathcal{H}$ sorted from closest to $\ell$ to furthest from $\ell$;
\STATE
Find an initial feasible solution, and solve the restricted version of the continuous relaxation of problem~\eqref{eqn:arc second-stage select} with the decision variables that are nonzero in the initial feasible solution;
\STATE
optimality$\_$verified $\leftarrow$ false;
\WHILE{not optimality$\_$verified}
    \STATE
    optimality$\_$verified $\leftarrow$ true;
    \STATE
    Use the optimal dual variables for the restricted version of the continuous relaxation of problem~\eqref{eqn:arc second-stage select} to compute $\check{c}_{t}(\ell)$ and $\check{h}_{t}(a,h',\ell)$;
	\FOR{$t \in \{1,..,T\}$}
	    \FOR{$a \in \mathcal{A}$}
            \FOR{$\ell \in \mathcal{L}$}
                \STATE
                negative$\_$found $\leftarrow$ false;
	            \FOR{$h' \in \mathcal{H}(\ell)$ (from closest to furthest) while not negative$\_$found}
                	\IF{objective value of column generation subproblem~\eqref{eqn:columngenerationsubproblem6} is negative}
                        \STATE
	                    negative$\_$found $\leftarrow$ true;
                        \STATE
                        optimality$\_$verified $\leftarrow$ false;
                        \STATE
                        Add variable $x_{t}(\check{c}_{t}(\ell),a,h',\ell,\check{h}_{t}(a,h',\ell))$ to the restricted version of the continuous relaxation of problem~\eqref{eqn:arc second-stage select};
                    \ENDIF
                \ENDFOR
            \ENDFOR
        \ENDFOR
    \ENDFOR
    \IF{not optimality$\_$verified}
        \STATE
        Solve the current restricted version of the continuous relaxation of problem~\eqref{eqn:arc second-stage select};
    \ENDIF
\ENDWHILE
\end{algorithmic}
\label{alg:columngeneration3}
\end{algorithm}

\subsection{Column Generation for Variables $x_{t}(c,a,\ell_{1},b,\ell,h)$}

The column generation subproblem to determine the variable $x_{t}(c,a,\ell_{1},b,\ell,h)$ with the smallest objective coefficient in the Lagrangian relaxation is
\begin{eqnarray}
\min & & \big\{f_{t}(c,a,\ell_{1},b,\ell,h) + \beta_{t}(a,b) \mathds{1}_{\{L(t,a,\ell_{1},b) = b\}} - \alpha_{t + \tau(t,c,a,\ell_{1},\ell,h)}(a,h) \nonumber \\
& & \hspace{5mm} {} + \psi_{t}(a,b,L(t,a,\ell_{1},b)) \mathds{1}_{\{L(t,a,\ell_{1},b) \in \mathcal{L}(t,a,b)\}} + \phi_{t}(c,\ell) + \theta_{t}(a,b,\ell_{1}) + \zeta_{t}(c,\ell) \; : \;
\label{eqn:columngenerationsubproblem1} \\
& & \hspace{10mm} t \in \{1,\ldots,T\}, c \in \mathcal{C}, a \in \mathcal{A}(c), b \in \mathcal{B}, \ell_{1} \in \mathcal{L}(t,a,b), \ell \in \mathcal{L}, h \in \mathcal{H}(c,\ell)\big\}, \nonumber
\end{eqnarray}
where the dual variable values are optimal dual values for the previous restricted problem.
Solving this column generation subproblem exactly can be time consuming, and is unnecessary for most iterations.
Next we show how to compute variables $x_t(c,a,\ell_{1},b,\ell,h)$ with negative reduced cost under the following simplifying assumptions:
\begin{enumerate}
\item[(A1)]
The dependence of the time $\tau(t,c,a,\ell_{1},\ell,h)$ on the emergency type~$c$, can be ignored.
Thus this time is denoted with $\tau(t,a,\ell_{1},\ell,h)$.
\item[(A2)]
The dependence of the cost $f_{t}(c,a,\ell_{1},b,\ell,h)$ on the ambulance station~$b$ or the emergency type~$c$, can be ignored.
Thus this cost is denoted with $f_{t}(a,\ell_{1},\ell,h)$.
\item[(A3)]
The set $\mathcal{A}(c)$ of ambulances to choose from does not depend on the emergency type~$c$.
Thus $\mathcal{A}(c) = \mathcal{A}$ for all~$c$.
\item[(A4)]
The set $\mathcal{H}(c,\ell)$ of hospitals to choose from does not depend on the emergency type~$c$.
Thus this set is denoted with $\mathcal{H}(\ell)$.
\end{enumerate}
Then the calculations can be streamlined as follows:
\begin{enumerate}
\item
For each time~$t$ and emergency location $\ell \in \mathcal{L}$, let
\[
c^*_{t}(\ell) \ \ \in \ \ \argmin\left\{\phi_{t}(c,\ell) + \zeta_{t}(c,\ell) \; : \; c \in \mathcal{C}\right\}
\]
denote the critical emergency type at time~$t$ and location~$\ell$.
\item
For each time~$t$, ambulance~$a \in \mathcal{A}$, and intermediate location $\ell_{1} \in \mathcal{L}$, let
\[
b^*_{t}(a,\ell_{1}) \ \ \in \ \ \argmin\left\{\beta_{t}(a,b) \mathds{1}_{\{L(t,a,\ell_{1},b) = b\}} + \psi_{t}(a,b,L(t,a,\ell_{1},b)) \mathds{1}_{\{L(t,a,\ell_{1},b) \in \mathcal{L}(t,a,b)\}} + \theta_{t}(a,b,\ell_{1}) \; : \; b \in \mathcal{B}\right\}
\]
denote the critical ambulance station at time~$t$ for ambulance~$a$ and intermediate location~$\ell_{1}$.
Observe that if $\left\{b \in \mathcal{B} \; : \; \ell_{1} \in \mathcal{L}(t,a,b)\right\} = \varnothing$, then there is no decision variable $x_{t}(c,a,\ell_{1},b,\ell,h)$ for such $\ell_{1}$.
\item
For each time~$t$, ambulance~$a \in \mathcal{A}$, emergency location $\ell \in \mathcal{L}$, and intermediate location $\ell_{1} \in \mathcal{L}$, let
\[
h^*_{t}(a,\ell_{1},\ell) \ \ \in \ \ \argmin\left\{f_{t}(a,\ell_{1},\ell,h) - \alpha_{t + \tau(t,a,\ell_{1},\ell,h)}(a,h) \; : \; h \in \mathcal{H}(\ell)\right\}
\]
denote the critical hospital at time~$t$ for ambulance~$a$, emergency location~$\ell$, and intermediate location~$\ell_{1}$.
\end{enumerate}
Then the column generation subproblem~(\ref{eqn:columngenerationsubproblem1}) reduces to
\begin{eqnarray}
\min & & \big\{f_{t}(a,\ell_{1},\ell,h^*_{t}(a,\ell_{1},\ell)) + \beta_{t}(a,b^*_{t}(a,\ell_{1})) \mathds{1}_{\{L(t,a,\ell_{1},b^*_{t}(a,\ell_{1})) = b^*_{t}(a,\ell_{1})\}} - \alpha_{t + \tau(t,a,\ell_{1},\ell,h^*_{t}(a,\ell_{1},\ell))}(a,h^*_{t}(a,\ell_{1},\ell)) \nonumber \\
& & \hspace*{5mm} {} + \psi_{t}(a,b^*_{t}(a,\ell_{1}),L(t,a,\ell_{1},b^*_{t}(a,\ell_{1}))) \mathds{1}_{\{L(t,a,\ell_{1},b^*_{t}(a,\ell_{1})) \in \mathcal{L}(t,a,b^*_{t}(a,\ell_{1}))\}} + \phi_{t}(c^*_{t}(\ell),\ell) \nonumber \\
& & \hspace*{5mm} {} + \theta_{t}(a,b^*_{t}(a,\ell_{1}),\ell_{1}) + \zeta_{t}(c^*_{t}(\ell),\ell) \; : \; t \in \{1,\ldots,T\}, a \in \mathcal{A}, \ell_{1} \in \mathcal{L}, \ell \in \mathcal{L}\big\}.
\label{eqn:columngenerationsubproblem2}
\end{eqnarray}
In most iterations, the column generation subproblem does not have to be solved to optimality.
Recall that in each column generation iteration, it is sufficient to find a variable $x_{t}(c,a,\ell_{1},b,\ell,h)$ with negative objective value in the column generation subproblem, or to verify that no variable $x_{t}(c,a,\ell_{1},b,\ell,h)$ has negative objective value.
Also, for any emergency location~$\ell$, an (available) ambulance at an intermediate location~$\ell_{1}$ that is close to $\ell$ is more attractive than an ambulance at an intermediate location~$\ell_{1}$ that is far from $\ell$ (this is also the intuition underlying the popular closest-available-ambulance dispatch heuristic).
These observations motivate Algorithm~\ref{alg:columngeneration1} below to solve problem~\eqref{eqn:columngenerationsubproblem2}.

\begin{algorithm}
\caption{Column generation algorithm for variables $x_{t}(c,a,\ell_{1},b,\ell,h)$}
\begin{algorithmic}[1]
\STATE
For each emergency location $\ell \in \mathcal{L}$, construct a list $\mathcal{L}(\ell)$ of intermediate locations $\ell_{1} \in \mathcal{L}$ sorted from closest to $\ell$ to furthest from $\ell$;
\STATE
Find an initial feasible solution, and solve the restricted version of
the continuous relaxation of problem~\eqref{eqn:arc second-stage select} with the decision variables that are nonzero in the initial feasible solution;
\STATE
optimality$\_$verified $\leftarrow$ false;
\WHILE{not optimality$\_$verified}
    \STATE
    optimality$\_$verified $\leftarrow$ true;
    \STATE 
    Use the optimal dual variables for the restricted version of the continuous relaxation of problem~\eqref{eqn:arc second-stage select} to compute $c^*_{t}(\ell)$, $b^*_{t}(a,\ell_{1})$, and $h^*_{t}(a,\ell_{1},\ell)$;
	\FOR{$t \in \{1,\ldots,T\}$}
	    \FOR{$a \in \mathcal{A}$}
            \FOR{$\ell \in \mathcal{L}$}
                \STATE
                negative$\_$found $\leftarrow$ false;
	            \FOR{$\ell_{1} \in \mathcal{L}(\ell)$ (from closest to furthest) while not negative$\_$found}
                	\IF{objective value of column generation subproblem~\eqref{eqn:columngenerationsubproblem2} is negative}
                        \STATE
	                    negative$\_$found $\leftarrow$ true;
                        \STATE
                        optimality$\_$verified $\leftarrow$ false;
                        \STATE
                        Add variable $x_{t}(c^*_{t}(\ell),a,\ell_{1},b^*_{t}(a,\ell_{1}),\ell,h^*_{t}(a,\ell_{1},\ell))$ to the restricted version of the continuous
                        relaxation of problem~\eqref{eqn:arc second-stage select};
                    \ENDIF
                \ENDFOR
            \ENDFOR
        \ENDFOR
    \ENDFOR
    \IF{not optimality$\_$verified}
        \STATE
        Solve the current restricted version of the continuous relaxation of problem~\eqref{eqn:arc second-stage select};
    \ENDIF
\ENDWHILE
\end{algorithmic}
\label{alg:columngeneration1}
\end{algorithm}


\newpage

Algorithm~\ref{alg:columngeneration4} is a single algorithm to find variables $x_{t}(c,a,b,\ell,h)$, $x_{t}(c,a,h',\ell,h)$, and $x_{t}(c,a,\ell_{1},b,\ell,h)$ with negative reduced cost for problem~\eqref{eqn:arc second-stage select}.

\begin{algorithm}[H]
\caption{Column generation algorithm for variables $x_{t}(c,a,b,\ell,h)$, $x_{t}(c,a,h',\ell,h)$, and $x_{t}(c,a,\ell_{1},b,\ell,h)$}
\begin{algorithmic}[1]
\STATE
For each emergency location $\ell \in \mathcal{L}$, construct a list $\mathcal{B}(\ell)$ of bases $b \in \mathcal{B}$ sorted from closest to $\ell$ to furthest from $\ell$;
\STATE
For each emergency location $\ell \in \mathcal{L}$, construct a list $\mathcal{H}(\ell)$ of hospitals $h' \in \mathcal{H}$ sorted from closest to $\ell$ to furthest from $\ell$;
\STATE
For each emergency location $\ell \in \mathcal{L}$, construct a list $\mathcal{L}(\ell)$ of intermediate locations $\ell_{1} \in \mathcal{L}$ sorted from closest to $\ell$ to furthest from $\ell$;
\STATE
Find an initial feasible solution, and solve the restricted version of the continuous relaxation of problem~\eqref{eqn:arc second-stage select} with the decision variables that are nonzero in the initial feasible solution;
\STATE
optimality$\_$verified $\leftarrow$ false;
\WHILE{not optimality$\_$verified}
    \STATE
    optimality$\_$verified $\leftarrow$ true;
    \STATE
    Use the optimal dual variables for the restricted version of the continuous relaxation of  problem~\eqref{eqn:arc second-stage select} to compute $\hat{c}_{t}(\ell) = \check{c}_{t}(\ell) = c^*_{t}(\ell)$, $\hat{h}_{t}(a,b,\ell)$, $\check{h}_{t}(a,h',\ell)$, $h^*_{t}(a,\ell_{1},\ell)$, and $b^*_{t}(a,\ell_{1})$;
	\FOR{$t \in \{1,\ldots,T\}$}
	    \FOR{$a \in \mathcal{A}$}
            \FOR{$\ell \in \mathcal{L}$}
                \STATE
                negative$\_$found $\leftarrow$ false;
                \FOR{$b \in \mathcal{B}(\ell)$ (from closest to furthest) while not negative$\_$found}
                	\IF{objective value of column generation subproblem~\eqref{eqn:columngenerationsubproblem4} is negative}
                        \STATE
	                    negative$\_$found $\leftarrow$ true;
                        \STATE
                        optimality$\_$verified $\leftarrow$ false;
                        \STATE
                        Add variable $x_{t}(\hat{c}_{t}(\ell),a,b,\ell,\hat{h}_{t}(a,b,\ell))$ to the restricted version of the continuous relaxation of problem~\eqref{eqn:arc second-stage select};
                    \ENDIF
                \ENDFOR
	            \algstore{cg_alg}
                \end{algorithmic}
                \label{alg:columngeneration4}
                \end{algorithm}
                \begin{algorithm}[H]
                \begin{algorithmic}[1]
                \algrestore{cg_alg}
                \FOR{$h' \in \mathcal{H}(\ell)$ (from closest to furthest) while not negative$\_$found}
                	\IF{objective value of column generation subproblem~\eqref{eqn:columngenerationsubproblem6} is negative}
                        \STATE
	                    negative$\_$found $\leftarrow$ true;
                        \STATE
                        optimality$\_$verified $\leftarrow$ false;
                        \STATE
                        Add variable $x_{t}(\check{c}_{t}(\ell),a,h',\ell,\check{h}_{t}(a,h',\ell))$ to the restricted version of the continuous relaxation of problem~\eqref{eqn:arc second-stage select};
                    \ENDIF
                \ENDFOR
	            \FOR{$\ell_{1} \in \mathcal{L}(\ell)$ (from closest to furthest) while not negative$\_$found}
                	\IF{objective value of column generation subproblem~\eqref{eqn:columngenerationsubproblem2} is negative}
                        \STATE
	                    negative$\_$found $\leftarrow$ true;
                        \STATE
                        optimality$\_$verified $\leftarrow$ false;
                        \STATE
                        Add variable $x_{t}(c^*_{t}(\ell),a,\ell_{1},b^*_{t}(a,\ell_{1}),\ell,h^*_{t}(a,\ell_{1},\ell))$ to the restricted version of the continuous relaxation of problem~\eqref{eqn:arc second-stage select};
                    \ENDIF
                \ENDFOR
            \ENDFOR
        \ENDFOR
    \ENDFOR
    \IF{not optimality$\_$verified}
        \STATE
        Solve the current restricted version of the continuous relaxation of problem~\eqref{eqn:arc second-stage select};
    \ENDIF
\ENDWHILE
\end{algorithmic}
\end{algorithm}


\end{APPENDIX}

\end{document}